\documentclass[12pt]{amsart}
\usepackage{amsmath}
\usepackage{amsfonts}
\usepackage{amssymb}

\def\ca{{\mathcal A}}
\def\cb{{\mathcal B}}
\def\cc{{\mathcal C}}
\def\cd{{\mathcal D}}

\def\ch{{\mathcal H}}

\def\ck{{\mathcal K}}

\def\cam{{\mathcal M}}
\def\cn{{\mathcal N}}

\def\cp{{\mathcal P}}
\def\cq{{\mathcal Q}}

\def\cS{{\mathcal S}}
\def\ct{{\mathcal T}}

\def\cx{{\mathcal X}}
\def\cy{{\mathcal Y}}
\def\cz{{\mathcal Z}}

\def\bc{{\mathbb C}}

\def\bn{{\mathbb N}}

\def\br{{\mathbb R}}

\def\bt{{\mathbb T}}

\def\eps{\varepsilon}

\def\k{\kappa}

\def\f{\varphi}

\def\Tr{\mathrm{Tr}}
\def\dist{\mathrm{dist}}
\def\dom{\mathrm{dom}}
\def\diam{\mathrm{diam}}

\newtheorem{Thm}{Theorem}[section]
\newtheorem{Cor}[Thm]{Corollary}
\newtheorem{Prop}[Thm]{Proposition}
\newtheorem{Lemma}[Thm]{Lemma}

\newtheorem{Dfn}[Thm]{Definition}
\newtheorem{exmp}[Thm]{Example}

\newtheorem{rem}[Thm]{Remark} 

\begin{document}
\title{{Extensions and Degenerations of Spectral Triples}}

 \author{Erik Christensen, Cristina Ivan}

\address{  
 Department of Mathematics, University of Copenhagen,   DK-2100\newline
\indent Copenhagen, Denmark}
 \email{echris@math.ku.dk}

 \address{   
Department of Mathematics, Leibniz University of Hannover,  30167\newline
\indent  Hannover, Germany }
 \email{ivan@math.uni-hannover.de}

 \date{\today}



 \begin{abstract}
   For a unital C*-algebra $\ca,$ which is equipped with a spectral
   triple $(A, H, D)$ and an extension $\ct$  of $\ca$ by the
   compacts, we construct a two parameter family of spectral triples
   $(A_t, K, D_{\alpha, \beta})$ associated to $\ct.$

Using Rieffel's notation of quantum Gromov-Hausdorff distance between
compact quantum metric spaces it is possible to define a metric on
this family of spectral triples, and we show that the distance
between a pair of  spectral triples varies
continuously with respect  the parameters.
It turns out that  a spectral
triple associated to the unitarization of the algebra of compact
operators is obtained under the  limit - in this metric - for  $(\alpha,
1) \to (0, 1),$ while the  
basic spectral triple $(\ca, H, D)$ is obtained from this
family under a sort of a dual limiting process for  $(1, \beta ) \to (1, 0).$  

We show that our constructions will provide families of spectral
triples for the unitarized compacts and for the Podle{\`s} sphere. In
the case of the compacts we investigate to which extent our proposed
spectral triple satisfies  Connes' 7 axioms for noncommutative geometry, \cite{Co3}.  
\end{abstract}

\maketitle

\section*{Introduction}
The so called Toeplitz algebra, say $\ct$, may be obtained in a number
of different ways. The most simple description of it is possibly as the
C*-algebra on the Hilbert space $\ell^2(\bn)$ generated by the
unilateral shift. A more profound description which relates to
analysis, can be obtained via the algebra, $\cc :=$
C$(\bt)$, of continuous functions on the unit circle. A function $f$
in this algebra is represented as a multiplication operator, $M_f$ on
the Hilbert space $H : = L^2(\bt)$ of square integrable functions.
This space has a subspace $H_+ $, which consists of those functions in
$H$ that have an analytic extension to the interior of the unit disk.
Let $P_+$ denote the orthogonal projection of $H $ onto $H_+$, then
the compression to $H_+$ of a multiplication operator $M_f$ for a
 continuous function $f$ on $\bt$ becomes the
Toeplitz operator $T_f : = P_+M_f|H_+,$ and these operators 
 form a subspace in
the Toeplitz algebra such that the Toeplitz algebra becomes the direct
sum of $\{T_f \, |\, f \in \cc \}$ and the algebra of compact
operators on $H_+.$ In this way the Toeplitz algebra becomes an
extension of $\cc$ by the compact operators. The mapping $ \cc \ni f \to
T_f$ relates to the differentiable structure on the circle 
 in the way that for the
ordinary differentiation on the circle with respect to arc length,
i. e.  
$D:= \frac{1}{i}\frac{d}{d \theta},$ we know that the space
$H_+$ is the closed linear span of the eigenvectors corresponding to
non negative eigenvalues for $D,  $   so there is a strong connection
between the differentiable structure on the circle and the
operator theoretical construction called {\em extension,} of $\cc$ by
the compacts.  In this article we will study this process from
a more general point of view. Our study is based on Connes' notion of
{\em a spectral triple } which is a way of expressing a differentiable
structure in the world of non-commutative *-algebras, \cite{Co2}.  

\begin{Dfn} \label{spectrip}
Let $\ca$ be a unital C*-algebra, $H$ a Hilbert space which carries a
faithful unital representation $\pi$ 
of $\ca$ and $D$ an unbounded self-adjoint
operator on $H$. For a dense  self-adjoint subalgebra $A$ of $\ca$
the  set $(A, H, D) $ is called a spectral triple associated to $\ca$  if
\begin{itemize}
\item[(i)] For all $a$ in $A$ the commutator $[\, D, \, \pi(a)\, ] $
  is bounded and densely defined.
\item[(ii)] the operator $(I+D^2)^{-1}$ is compact.
\end{itemize} 
\end{Dfn}

In this article our starting point is a spectral triple associated to
a C*-algebra $\ca$  and we want to
study some of the possibilities for constructing spectral triples
associated to an extension of $\ca$ by the algebra of compact
operators on an infinite dimensional
 Hilbert space. Our fundamental example of a spectral
triple is the
one coming from the unit circle, as described above.
This particular example was investigated by Connes and Moscovici in
\cite{CM}, where  they constructed a spectral triple associated to the
the Toeplitz algebra for each natural number $n$  in 
the following way. Let $S$ denote the unilateral shift on $H_+, $
i. e. for the canonical basis for $H_+$ we have $Se_k = e_{k+1},$ and
let $D_+ = -i \frac{d}{d \theta}| H_+,$ which means that $D_+$ is the
positive self-adjoint operator on $H_+$ which satisfies $D_+e_k = k
e_k.$  Then
the Hilbert space $K$ of the spectral triple for the Toeplitz algebra
is defined as $K \, = \, H_+ \oplus
H_+,$ 
and the Dirac operator $D_n$ is defined via the matrix form
\begin{displaymath} D_n \, :=  \, \begin{pmatrix} 0 & D_+ S^n
    \\ (D_+ S^n)^* & 0 \end{pmatrix}. \end{displaymath}
In the construction we present in this paper we look at a spectral
triple $(A, H, D)$ associated to a C*-algebra $\ca$ and an orthogonal
projection $P$ onto a subspace of $H$ such that $P$ commutes with $D,$
and for each operator $a$ from $\ca$ the commutator $[P, a] $
is compact. All of this set-up  
 is analogous to the classical spectral triple for the
circle algebra, but there is, in general, no counterpart to the
unilateral shift. This means that we have to modify the construction by
Connes and Moscovici in order to construct a spectral triple
associated to the C*-algebra generated by the operators $\{Pa|PH \,
\big| \, a \in \ca\,\}$ and the compact operators on $PH.$ A C*-algebra,
obtained this way, is called an extension of $\ca$ by the compacts, and
one of the  problems we try to solve in this article is to
find ways to {\em extend } a spectral triple associated to a
C*-algebra to a spectral triple for an extension of that algebra by
the compacts. A more general question of this sort has been studied by
Chakraborty in \cite{Ch}. In that paper he studies compact quantum
metric spaces as introduced by Rieffel \cite{RiStSp}, and he investigates
the possibilities to generate the structure of a compact quantum metric space
associated to an extension of a C*-algebra which is associated to the
given  compact quantum metric space. 
 In Chakraborty's article he 
 studies a short exact sequence of C*-algebras of the type 
\begin{displaymath}   
0 \rightarrow \ck \otimes \ca \rightarrow \ca_1 \rightarrow \ca_2 \rightarrow 0,
\end{displaymath}
for which the last homomorphism has a positive splitting $\sigma : 
\ca_2 \to \ca_1,$  and he shows that if there is a compact quantum
metric space associated to both $\ca$  and   $\ca_2$, then there
exist several  compact quantum
metric spaces associated to the C*-algebra $\ca_1.$ This set up is more
general than ours from the point of view of  possible extensions, but our
concern is spectral triples rather than the construction of compact
quantum metric spaces.  Chakraborty
offers several applications of his construction to known examples, such
as the Podle{\`s} sphere. We show that our
construction can also be applied to generate spectral triples for 
this example and also for the algebra of
compact operators on a separable Hilbert space.

We will, through the entire article, suppose that $(A, H, D)$ is a
spectral triple associated to a unital C*-algebra $\ca,$
 which is a subalgebra of $B(H).$ As in the book
\cite{HR} Definition 2.7.7 and Chapter 5, we will study extensions of
Toeplitz type. This means that we are interested in an
orthogonal projection  $P$ in $B(H)$ which commutes with $\ca$ modulo
the compact operators.  One can then define 
a C*-algebra $\cb$ on $PH$ as the  C*-algebra  generated by the space
of operators  $ \{
Pa|PH\, \,  \big| \, a \in \ca \,\}$ in $B(PH).$ For each operator $a$ in $\ca$ we
have that $(I-P)a|PH$ is compact, so operators of the form
$Pa^*(I-P)a|PH$ are compact and - unless $P$ commutes with $\ca$ - the
algebra $\cb$ will contain non trivial compact operators. We will let
$C(PH)$ denote the algebra of compact operators on $PH.$ 
 In the classical case of the Toeplitz algebra for the circle 
we actually have $C(PH) = \cb \cap C(PH). $ In
any case, independently of what the C*-algebra $\cb \cap C(PH) $ might
be we will consider the C*-algebra $\ct$ which is defined as the sum $\ct := \cb + C(PH).$

Let $\cq(PH)$ denote the Calkin algebra $B(PH) / (CPH)$ and let $\k$
denote the quotient homomorphism,  then we can  define a 
 homomorphism $\f$  of
$\ca$ into $\cq(PH)$  by $\f(a) := \k(Pa|PH).$
The extensions we will consider are those obtained via  the
construction described above which also have the property that the 
homomorphism $\f$  is faithful on $\ca$.
We may then  define a homomorphism of $\ct$  onto $\ca$ as $\f^{-1}\circ\k,$ 
and we will say that a projection $P$ in
$B(H)$ which satisfies all the properties discussed here is  
of Toeplitz type. We will not study all such projections, but
restrict our investigations to projections  of Toeplitz type such that  $P$ 
commutes with the Dirac operator $D$ and satisfies
 the following regularity property 

\begin{equation} \label{regular} 
\forall a \in A \quad [PD,a] \quad
\text{ is bounded and densely defined.}
\end{equation} 
\noindent
Under these assumptions we will say that the quadruple $\big((A,H,D),
P\big)$ is of Toeplitz type. 

 We would like to
remind the reader that for any spectral triple - of infinite dimensions
- like $(A, H, D)$ 
the spectral projection $P_+$ for $D$ corresponding to the
interval $[0, \infty[$ is a natural candidate for $P.$ 
 This follows from the well known fact that the
symmetry   $2P_+ -I$   generates a bounded Fredholm module 
\cite{BJ}. In the case $P=P_+$ the regularity condition amounts to the
assumption that for any $a$ from $A,$ the commutator $[|D|, a ] $ is
bounded, too. 

We can now describe the general construction,
which we will study here. So for a  quadruple $\big((A,H,D),
P\big)$  of Toeplitz type  we define the
C*-algebra $\ct,$ as above, a  Hilbert
space $K := PH \oplus H $  and a representation $\pi$ of $\ct$ on $K$
given by the matrix form
  \begin{displaymath} 
\pi(t) \,:= \begin{pmatrix} t & 0 \\ 0 & \f^{-1}( \k(t))\end{pmatrix}.
\end{displaymath}

The Hilbert space $K = PH \oplus H$   decomposes as   $K = PH
\oplus PH \oplus (I-P)H $ and we can see that the first two summands
are exactly analogous to the ones appearing in the Connes-Moscovici
construction. We have tried to follow their idea, but our analysis indicates
that we can not in general give up the information which is encoded in
the $(I-P)H$ part of the structures, so we will consider $K = PH
\oplus H$ rather
than $PH \oplus PH.$ On the other hand this opens the possibility to
play on the two parts with {\em different weights} as the introduction
of parameters in our proposal for a Dirac operator shows.
Since $D$ is supposed to commute
with $P,$ the regularity conditions imposed make it possible to 
 define a family of Dirac operators on $K$ in the following way. 
For positive reals $\alpha,
\beta $ such that $\alpha \beta \leq 1 $ we define 
an unbounded self-adjoint
operator $D_{\alpha, \beta} $ on $K$ via its matrix:
\begin{displaymath}
D_{\alpha, \beta} \, := \, \begin{pmatrix} 0 & \beta D|PH & 0\\
\beta D|PH   & \frac{1}{\alpha} D|PH  & 0 \\ 0 & 0 & \frac{1}{\alpha}
D|(I-P)H \end{pmatrix}
\end{displaymath}

 The reason for having the
parameter $\alpha $ appearing in the form $1/\alpha$ is mainly {\em
  aesthetical.  } For instance, the formula giving a  distance
estimate between the non-commutative spaces obtained for two pairs
of parameters, say $(\alpha, \beta) $ and $(\gamma, \delta) $ 
becomes by Theorem \ref{distThm}

\begin{displaymath}
\big( \max \big\{\frac{\alpha \beta }{\gamma  \delta }, \frac{\gamma  \delta }{\alpha \beta } \big\} -1 + \big|1-
\frac{\beta}{\delta }\big|\big)\diam_{\alpha, \beta },
\end{displaymath}
and in this formula the product $\alpha \beta$ fits in as a parameter.

The reason why the parameters are supposed to satisfy the inequality
$\alpha \beta \leq 1 $ comes originally from the classical Toeplitz
case, where it is quite easy to analyze the situation in details.
It turns out that for this
example and a pair of parameters $(\alpha, \beta) $ for which  $\alpha \beta
> 1 $ all the aspects of the non-commutative space associated to these
parameters is already contained in the space given by the parameters
$(1/\beta, \beta).$ The general case does not work exactly in the same
way, but the Remark \ref{leq1}, explains why we think nothing essential
is lost, if we just stick to the region in the parameter space where
$\alpha \beta \leq 1.$ Still another argument for the choice of
parameter space is that it turns out that the limiting process $\alpha
\to 0 $ behaves uniformly nice on the set of parameters where $\beta
\geq \beta_0.$ For the convergence $(\alpha, 1) \to (0,1) $ we find
that it induces a convergence with respect to the quantum
Gromov-Hausdorff metric of the compact quantum metric spaces
associated to $\ct$ for the parameters $(\alpha, 1)$ onto that of the
unitarized algebra $\tilde \cc.$ Our intuitive description of
this phenomenon is as follows. We think that the spectral triple 
acts like a
microscope where we have a fixed screen to watch, but we are allowed
to change the magnification. The parameter $\alpha$ is a measure of
the actual size, say in meters,  of the objects we can watch on our
screen, and 
 then  $1/\alpha$ is the
magnification factor. When $\alpha$ decreases to zero, we loose the
sight of the big picture and can only see tiny details of very small
things. In the end - when $\alpha = 0\, -$ 
our mathematical construction can only see the compacts, and they are
considered to be the {\em infinitesimals } in Connes' dictionary
\cite{Co3}.
The precise mathematical content of this {\em story } is contained in
Theorem \ref{ato0}. 

The limit $(1,\beta) \to (1,0)$ is quite easy to
understand if you take a look at the definition of $D_{\alpha, \beta}
$ just above, and you can see that you get the basic
spectral triple $(A, H, D)$ back, but now in a degenerated representation. 
This is not a limit with respect to the quantum Gromov-Hausdorff
metric on the associated non commutative spaces, but rather a sort
of degenerated limit where the compacts i. e. the infinitesimals
become invisible. We can provide a simple 2-dimensional model of the
picture we try to present. Look at the unit square $[0,1]^2$ in
$\br^2$ and equip it with the metric $d_{\alpha, \beta}\big((x,y),
(s,t)\big):= \alpha|x-s| +  1/\beta|y - t|,$ then the limit $(\alpha, 1)
\to (0,1)$ gives the unit interval $ \{(0,y)\,|\, 0 \leq y \leq 1 \} $
with its standard metric 
as the limit in the Gromov-Hausdorff metric.
 For $(1, \beta) \to (1,0)$
 there is no limit of this sort but we get - pointwise -  a degenerate metric
 $d_{1,0}$  on the unit square as a limit. This degenerate metric is
 given by the formula  
\begin{equation} \label{d2}
d_{1,0 }\big((x,y),(s,t)\big)\, = \, \begin{cases} |x-s| \, \text{ if } t
  = y
\\ \infty \quad \quad \text{ if } t \neq y 
\end{cases}.\end{equation}

In the first version of the article we considered this limiting
process as a sort of deformation, but we have been told by several
people that we do not deform a product, so the wording is wrong. By
looking into the literature on metric spaces associated to Riemannian
structures we have found that the limiting processes we are watching
may be considered as degenerations of metric spaces, so the title has
been changed accordingly. Usually this sort of degeneration of
Riemannian structures is studied  under some assumptions on
boundedness of curvature during the process \cite{Be,BBI,Fu}, but this
last aspect of degeneration of metric structures does not apply
to  our results, at least for the time being.

On the other hand we still think that the limiting process $(1, \beta)
\to (1,0) $ - in the case where the algebra $\ca$ is commutative -
offers a way of describing a {\em passage } from a non-commutative
compact metric space into a commutative compact metric space. The
equation (\ref{d2}) indicates  that a better description of the
degeneration occurring while $\beta \to 0$ might possibly be;  
{\em a passage from a non-commutative space to an infinite collection
  of disjoint identical copies of the same commutative space.}

 After we had posted the first pre-print version of this article on
 the arXiv, we were informed by Hanfeng Li, that our constructions can
 work in the settings of David Kerr's, \cite{Ke}, and his own,
 \cite{Li}, and the one from their joint work, \cite{KL}, where the
 quantum Gromov-Hausdorff metric is extended in a way which is based on the
 operator space structure of the given algebra. The 
 introduction of  the state spaces of the tensor products of the given
 algebra by the algebras of $n \times n$ complex matrices into the
 definition of a compact quantum space, is to be able to describe
 certain aspects of order in more
 details. We have chosen not to expand the present paper and hope that
 these results of Li's some day will find a suitable place to be
 presented.  

Near the end we give a couple of examples and show that
 our method creates an abundance of spectral triples for the
 unitarized compact operators on a separable infinite dimensional
 Hilbert space. Then we show that the method, when applied to the unit
 circle and the classical differential operator
 $\frac{1}{i}\frac{d}{d\theta},$ gives a spectral triple associated to
 the classical Toeplitz C*-algebra. Based on this spectral triple we
 can then by a slight modification of our method obtain a spectral
 triple for the Podle{\`s} sphere. Unfortunately we can see no
 relations between our constructions and the ones presented in
 \cite{DDLW} and \cite{DS}.

At the very end we present some small comments on the relations
between the constructions in this  article  to  the concepts of even and
odd spectral triples and  to  analytic K-homology as described by Higson
and Roe in their book \cite{HR}. It may be that further assumptions
or conditions on the starting spectral triple may be used to give a
basis for a more detailed study of such relations. In this paper we
have been focusing on the quite general degeneration aspects of the extended
spectral triples.

We are most thankful to the referee who has pointed out some problems
in the first version of the article, and he has also suggested several
possibilities for improvements. Among the questions he asked is the
question on how the spectral triples constructed here relate to the 7
axioms for non-commutative geometry which Connes lists in the article
\cite{Co3}. To answer this question we have studied this question for  our examples
involoving the Toeplitz algebra and the unitarized compacts. The
Toeplitz
 case seems not promising at all with respect to this investigation,
 so we have not included any comments on this aspect for the Toeplitz
 algebra. For the compacts we do check all the axioms, and we
show that we can meet some, whereas for others we can not decide, but for the
so-called reality axiom  we get one of the  signs wrong.

\section{A family of spectral triples associated to an
  extension}

We will keep a C*-algebra $\ca$  with an associated spectral triple
\newline 
$(A, H, D)$ fixed during the
whole article and moreover
suppose that $\ca$ is a concrete C*-algebra acting on the Hilbert
space $H.$ As stated in the introduction, we will assume that we have a
projection $P$ in $B(H)$ of Toeplitz type and study an extension
$\ct$ of $\ca$ by the algebra  $C(PH)$ of compact
operators  on $PH.$ 
 It should be remarked that {\em  we do
not assume that the algebra of compact operators $C(PH)$ is contained
in the C*-algebra $\cb$ generated 
by operators of the form $Pa|PH,$ }  and in particular we can also study
the situation where $P$ commutes with $\ca.$ 
We will collect the definitions from the introduction in a formal
definition. 

\begin{Dfn}  \label{Ttype}
Let $\ca$ be a unital C*-algebra on a Hilbert space $H$ and let 
\newline $(A, H, D) $ be a spectral triple associated to $\ca$. A
projection $P$ in $B(H)$ is said to be of Toeplitz type for $(A,
H, D)$ if 
\begin{itemize} 
\item[(i)] The projection $P$ commutes with $D.$
\item[(ii)] The projection $P$ commutes modulo the compacts with
  $\ca.$
\item[(iii)] The homomorphism $\f$ of $\ca$ to the Calkin algebra
  $\cq(PH),$ defined by $\f(a) := Pa|PH + C(PH),$ is faithful.
\end{itemize} 
For such a triple and a projection $P$ of Toeplitz type we define the
Toeplitz extension $\ct$ of $\ca$ by $C(PH)$ as the C*-algebra generated by
\begin{displaymath} \{ PA|PH \, \big| \, A \in \ca \} \cup C(PH)
\end{displaymath} 
\end{Dfn}

At the end of the paper we construct an example which will give a
spectral triple for the Podles' sphere. That example is based on a
slight variation of the construction presented in this paper, and it 
 suggests that it might be possible to study extensions of $\ca$  by a
 sub C*-algebra of $C(PH)$ instead.  A
generalization of our construction to cover cases like this seems
possible, but also quite demanding with respect to extra details, so
we have chosen only to consider extensions by all of $C(PH),$
 and then just present the other point of view in connection
with the example for the Podles' sphere.

\medskip

Given a projection $P$ of Toeplitz type for $(A, H, D ),$   
we assume that $P$ commutes with $D$ and by this we mean that $P$
commutes with all the spectral projections of $D.$ From this it
follows that the unitary $S:= P - (I-P)$ also will commute with $D$
and $S$ will  map the  domain of definition for $D$ 
 onto itself, ( \cite{Pe},
Proposition 5.3.18.) Hence the domain of
definition for $D$ splits into a direct sum of its intersections with $PH$
and $(I-P)H$ respectively. We will need  that the commutators from the
spectral triple respect this decomposition too, and this is the basis
for the following definition. In the classical case where $P= P_+$
this means that we will not only demand that commutators $[D, a] $ are
bounded and densely defined for $ a $ in $A,$ but we want both $[D,a]
\text{ and } [|D|, a]$ to be bounded and defined on a common dense
domain.

\begin{Dfn} \label{Tspt} 
A quadruple  $\big( (A, H, D),P \big) $ where $P$ is a projection of
Toeplitz type for   
$(A, H, D), $ is said to be of Toeplitz type if:

\begin{itemize}  
\item[(i)] For any  $a$ in $A$, the   commutators 
  $[PD,a]$ and $[(I-P)D,a] $ are bounded and densely defined and their
  common domain of definition
  contains two subspaces \begin{displaymath}\dom ([D,a]) \cap PH
    \text{ and }
  \dom([D,a])\cap (I-P)H \end{displaymath}
 which are dense in $PH$
  and $(I-P)H$ respectively.
\item[(ii)] The operator $D_{P} := D| PH$ has trivial kernel.
  \end{itemize} 
\end{Dfn}

The properties  in the definition above seem natural in the  setting
for a classical Toeplitz algebra, except for the last one. On the
other hand that one does not really matter.
 Let namely $N$ denote the orthogonal projection onto the kernel
of $D,$ then $N$ is of finite rank, and since it is a spectral
projection for $D,$ it commutes with $P$ and we can replace $P$ by
$P- PN, $ without disturbing any properties of the extension we are
studying. 
The first condition has been imposed 
in order to be able
to look at commutators of the form $[PD, a]|PH$ and their relatives
with restrictions to $(I-P)H$ and/or $PD$ replaced by $(I-P)D.$ The
conditions are made such that the lemma below holds.
To keep the notational
problems at a minimum we introduce the conventions that
\begin{align*}
H_p &:= PH, \,\, H_q:= (I-P)H, \\ \, P_p &:= P, \,\quad P_q:= (I-P),\\
\, D_p& := DP,\, \,  D_q:= D(I-P).\end{align*}

\begin{Lemma} \label{closcom}
For any  $a$ in $A$ and  any combination of the symbols  $s,t,r$ in the set
$\{p,q \}$ 
\begin{align}
\text{The closure of }\big(P_s [D,a]|H_t\big)\,& = \,P_s\text{the closure of }\big(
[D,a]\big) |H_t  \label{Dcom} \\ 
\text{The closure of }\big(P_s [D_r,a]|H_t\big)\,& = \,P_s\text{the closure of }\big(
[D_r,a]\big) |H_t  \label{|D|com} 
\end{align}
\end{Lemma}
\begin{proof}
We will not prove all these statements but restrict ourselves to the
relation (\ref{Dcom}) in the situation where $s = p$ and $t = q.$ 
The closure of the commutator $[D,a]$  
is bounded and we will denote its closure by $\delta(a).$ It is
immediate that as operators we have the inclusion  $P[D,a]|(I-P)H
\subseteq P\delta(a)|(I-P)H, $ and in order to show the statement of the
lemma it is sufficient to show that $P[D,a]|(I-P)H$ is densely defined,
but this is fulfilled by the condition (i) in Definition \ref{Tspt}. 
We now claim that we can perform exactly the same computations with
respect to any other combination of the symbols $\{p, q\},$ and then
obtain the lemma. 
\end{proof}

\noindent
The effect of the lemma is that we may decompose the commutator
$[D,a]$ into its matrix parts with respect to the decomposition of $H=
H_p  \oplus H_q,$ such that each of the 4 the matrix entries of the
closure is the closure of the corresponding operator-theoretical
matrix entry. From this follows the lemma just below:

\begin{Lemma} \label{offdiag}
For any $a$ in $A$ the operators $DP_paP_q$ and $DP_qaP_p$ are 
  bounded and everywhere defined. 
\end{Lemma}   

\begin{proof}
We remind you that a product of operators of the form $CB$ where $C$
is closed and $B$ is bounded is automatically closed, so if it is
bounded it must be everywhere defined.
\end{proof}

\medskip

We will now define various maps and a
spectral triple associated to $\ct.$
Before we give the
definition we would like to mention that its first  item is
legal, due to a general result on
ideals in C*-algebras that we recall here (\cite{Da}, Corollary 1.5.6). 
\begin{Prop}
Suppose that $\mathcal{I}$ is a two sided closed 
 ideal of a C*-algebra $\mathcal{A}$, and that $\mathcal{B}$ is a
 sub C*-algebra of $\mathcal{A}$. Then $\mathcal{B}+\mathcal{I}$ 
is a C*-algebra and 
\begin{displaymath}
\mathcal{B}/(\mathcal{B}\cap \mathcal{I})\simeq (\mathcal{B}+\mathcal{I})/ \mathcal{I}
\end{displaymath}
is a *-isomorphisms. 
\end{Prop}

Let now $\k$ denote the quotient mapping of $B(PH)$ onto the Calkin
algebra $\cq(PH),$ then the proposition above has the following
corollary as a consequence.

\begin{Cor}
For any  quadruple $\big( (A, H, D),P \big) $ of Toeplitz type 
with associated 
Toeplitz
extension $\ct, $ the images  $\k(\ct)$ and $ \varphi(\ca)$ in the Calkin
algebra $\cq(PH)$ agree and $\k(\ct)$ is
isomorphic to $\ca.$
\end{Cor}

\begin{Dfn} \label{Tmap}
Let $\big((A,H,D), P \big)$ be a quadruple of Toeplitz type associated
to a C*-algebra $\ca.$ For the induced Toeplitz extension $\ct $ of
$\ca$ we define:
\begin{itemize} 
 
\item[(i)] A representation $\rho : \ct \to B(H)$ by $\rho(t) :=
  \f^{-1}(\k(t)).$ 
\item[(ii)] A completely positive unital and injective mapping 
$T: \ca \to \ct$ by $T(a):= Pa|PH.$ 
\item[(iii)] A projection $\Theta$ of  $\ct$ onto $C(PH)$ by 
$\Theta(t) := t - T(\rho (t)). $
\item[(iv)] For any $ x $ in $B(H)$ and any combination of the symbols
    $s,r \in \{p, q\} $ we define $x_{sr}$  in $B(H_r, H_s)$ by 
$x_{sr} := P_sx|H_r.$ 
\end{itemize}
\end{Dfn}
\noindent
It should be noted that for an $a$ from $\ca,$ we have $T(a) = a_{pp}.$ 

Given a situation as above, we will then define a representation $\pi$
of $\ct$ on a Hilbert space $K$ and a family of unbounded self-adjoint
operators $D_{\alpha, \beta} $ on $K,$ but it is not immediate that we
will get spectral triples this way, so we start by defining the
ingredients separately and study some of their properties.

\begin{Dfn}  \label{ST}
\indent Let $\big((A,H,D), P \big)$ be a quadruple of Toeplitz type associated
to the C*-algebra $\ca$  and let
$\ct$ denote the induced Toeplitz algebra on the space $PH.$ To this
quadruple is associated:
\begin{itemize} 
\item[(i)] A dense self-adjoint subalgebra $A_c$ of $C(PH)$ defined
  by 
\begin{displaymath} A_c \, := \, \{ k \in C(PH)\, | \,  D_pk \text{
    and } kD_p \text{ are bounded }\, \} \end{displaymath}
\item[(ii)] A dense self-adjoint subalgebra $A_t$ of $\ct$ defined
  by
\begin{displaymath}  A_t \, := \, \{ T(a) + k \, |\, a \in A,\, \, k
  \in A_c \, \} \end{displaymath}
\item[(iii)] A Hilbert
space $K$ defined as the sum 
\begin{displaymath}K\, := \, PH \oplus H \, =
  \, H_p \oplus H_p  \oplus H_q \end{displaymath}

\item[(iv)] 
  A representation $\pi$ of  $\, \ct $ on $K$  defined by   
\begin{displaymath}
\forall t \in \ct: \quad \pi(t) \, := \, \begin{pmatrix} t &  0 \\ 0 &
  \rho(t) \end{pmatrix}.
\end{displaymath}
\item[(v)] 
For  positive reals $\alpha,  \beta  $  a  self-adjoint
operator $D_{\alpha,\beta }$  is defined  on $K $ via its 
 matrix, which, with respect to the
decomposition  $K \, = \,  H_p \oplus H_p \oplus H_q,$ is given by   
\begin{displaymath} 
D_{\alpha,\beta} \, := \, \begin{pmatrix}0 & \beta D_p & 0 \\ \beta D_p
  & \frac{1}{\alpha} D_p & 0 \\ 0 &
  0 &  \frac{1}{\alpha} D_q \end{pmatrix}.
\end{displaymath}
\end{itemize}
\end{Dfn}
\noindent
It may not be obvious that the linear space $A_t$ is an algebra, but
it follows from Lemma \ref{offdiag}.

We will show that for each pair $(\alpha, \beta)$ we will get a
spectral triple for the Toeplitz extension $\ct$  of $\ca,$ induced by the
projection $P.$ This will be an odd spectral triple and it is possible
- via a standard trick -  to
obtain  an even triple instead. But from the point of view we are
studying here, namely the variation of the compact quantum 
 metric spaces with respect
to the parameters $\alpha $ and $\beta $ we do not get any changes if
the investigation is performed with the odd spectral triple described
above or an even one.
 
The properties of a quadruple of  Toeplitz type
 now come into play and it helps
us to split a commutator $[D_{\alpha, \beta}, \pi(t)]$ for a $t$ in
 $A_t$ into its matrix parts.

\begin{Lemma} \label{ComT}
For $a$  in $A$ and $k$ in $A_c$ and positive reals
$\alpha, \beta$  with $\alpha \beta \leq 1 $ the commutator 
$[D_{\alpha, \beta}, \pi((T(a)+k) ]$  is bounded. For each matrix
part of the closure of 
this commutator, with respect to the decomposition $K = PH
\oplus PH \oplus (I-P)H$, the element is the closure of the
corresponding matrix part of the algebraic commutator. 
\end{Lemma}

\begin{proof}
We will do the  computations where they are defined purely
algebraically, then show that each matrix entry is bounded and densely
defined and then conclude that the closure of the commutator 
 is the closure of the
operator composed of the matrix entries. The reason why this is
possible is the regularity assumptions and the Lemma \ref{closcom}.

 \begin{equation} \label{TComm}  
 \end{equation}
\begin{align*} & \, [D_{\alpha, \beta}, \pi(T(a)+k)] \, \\ 
& \, =  \left [
   \begin{pmatrix} 0 & \beta D_p  & 0 \\ \beta D_p  & 
  \frac{1}{\alpha} D_p & 0 \\ 0 & 0 &  \frac{1}{\alpha} D_q
\end{pmatrix}, \, 
\begin{pmatrix} T(a)+k & 0 & 0 \\ 0 & T(a) & a_{pq}\\
  0 & a_{qp} & a_{qq}  \end{pmatrix} \right ]   \\ 
& \, =  \beta \begin{pmatrix} 0 & [ D_p, T(a)] - kD_p & D_pa_{pq} 
   \\ [D_p,T(a)] +D_pk & 0 & 0 \\ -a_{qp}D_p & 0 & 0 \end{pmatrix} \\
 & +
 \frac{1}{\alpha} \begin{pmatrix} 0 & \begin{matrix} 0 & 0  \end{matrix} 
\\ \begin{matrix} 0 \\ 0 \end{matrix} &  [D, a] \end{pmatrix}. 
\end{align*}  
 The lemma follows. 

\end{proof}

\begin{rem} \label{leq1}
The idea in the setup of the commutator \begin{displaymath}[D_{\alpha,
    \beta}, \pi(T(a) + k)] \end{displaymath}
 is that it shall reflect both of the
  given commutators $[D_p, k] $ and $[D, a]$ in such a way that a
  variation of the parameters $\alpha $ and $\beta$ will reveal
  information on each of these parts separately. Since
  \begin{align*} [D_p, T(a)] & = P[D, \pi(a)| PH  \text{ we get }
\\ \| [D_p, T(a)] \| & \leq  \|[D, \pi(a)]\| \text{ and then for }
\alpha \beta \leq 1 \\  
\beta \| [D_p, T(a)] \| & \leq \frac{1}{\alpha} \| [D_p, T(a)] \|
\leq \frac{1}{\alpha}  \|[D, \pi(a)]\|
\end{align*}
Hence for $\alpha \beta \leq 1 $ the term $\beta [D_p, T(a)]$ will not
be of significance and then we see from (\ref{TComm}) that in this
case we will have  \begin{align*}
\frac{1}{\alpha}  \|[D, \pi(a)]\| & \leq 1 \text{ and }  \max \beta
\{\|D_pk\|, \|kD_p\| \} \leq 2 \text{ if } \\  \|[D_{\alpha,
  \beta}, \pi(T(a) + k)]\| & \leq 1. 
\end{align*}  
We will then impose the condition $\alpha
\beta \leq 1 $ in all of our future statements, and 
 this also fits nicely with the results of Theorem   
 \ref{distThm} which indicate that the product $\alpha \beta $ is a
 relevant parameter. 
\end{rem}

\begin{Prop}  \label{SPT}
For any pair of positive real numbers $\alpha,  \beta $ such that
$\alpha \beta \leq 1 $  the
tuple 
\begin{displaymath}
(A_t , K, D_{\alpha, \beta} ) 
\end{displaymath}
is a spectral triple associated to the C*-algebra $\ct.$
This extended   spectral triple is $s-$summable if and only if the
given one is $s$-summable.
\end{Prop}

\begin{proof} 

Having the Lemma \ref{ComT} we just have to prove that each $D_{\alpha,
  \beta} $ has compact resolvents, but that follows immediately from
the definition of  $D_{\alpha,
  \beta }. $  Since
$P$ commutes with the spectral projections for $D,$ each eigenspace
$H_{\lambda_i} $ for $D$ decomposes as an orthogonal sum   $PH_{\lambda_i}
\oplus (I-P)H_{\lambda_i} $ and we can find an orthonormal basis for
$PH,$ say $(e_i),$ consisting of eigenvectors for $D_p,$ plus an
orthonormal basis for $(I-P)H,$  say $(f_j),$ consisting of
eigenvectors for $D_q.$ If $e_i$ is an eigenvector corresponding to
the eigenvalue $\lambda_i, $  the operator $D_{\alpha, \beta} $
will have an invariant 2 dimensional subspace of the form $\{(z e_i, w
 e_i, 0 ) \, |\, z, w \in \bc \, \}$  in the decomposition of $K.$ The
 eigenvalues of $D_{\alpha, \beta}$ on this space are determined
 by the $ 2 \times 2 $ matrix 
\begin{displaymath} 
M(\alpha, \beta) \, := \, \begin{pmatrix} 0 & \beta \\
\beta & \frac{1}{\alpha} \end{pmatrix}  
\end{displaymath}
\noindent 
such that the eigenvalues become $\lambda_i $ times the eigenvalues of
$M(\alpha, \beta).$  
For an eigenvector $f_j$ for $D_q$ corresponding to an eigenvalue
$\mu_j$  this vector becomes an eigenvector for $D_{\alpha, \beta}$
corresponding to the eigenvalue $\mu_j / \alpha .$ Let now $s$ denote a
positive real and we see that we get the equality  below

\begin{displaymath}
\Tr(|D_{\alpha, \beta} |^{-s}) \, = \Tr( |M(\alpha,
\beta)|^{-s})\, \Tr(|D_p|^{-s} )\,  +  \, \alpha^s \Tr(|D_q|^{-s})
\end{displaymath}
and the proposition follows.
\end{proof}

While we are at such matrix computations we remind you that for
positive real numbers $\alpha, \beta, \gamma, \delta, $ Hilbert spaces
$L, M, N$  and bounded operators $v \in B(M,L), \, x \in B(L,M), \, 
y \in B(M), \, z \in B(N),$ 
 we can obtain the identities below with respect to some  operator matrices
 on  Hilbert sum $L
 \oplus M \oplus N.$ 

\begin{align} \label{abcd}
&\begin{pmatrix} 0 & \beta v & 0 \\ \beta x & \frac{1}{\alpha}y & 0 \\
0 & 0 & \frac{1}{\alpha} z 
\end{pmatrix} \, = \\ & 
   \begin{pmatrix} \sqrt{\frac{\alpha}{\gamma}}
      (\frac{\beta}{ \delta} ) I
    & 0 & 0 \\ 0 & \sqrt{\frac{\gamma}{\alpha}}I & 0 \\ 0 & 0 
& \sqrt{\frac{\gamma}{\alpha}}I  
  \end{pmatrix} \begin{pmatrix} 0 & \delta v & 0 \\
  \delta x &  \frac{1}{\gamma} y & 0 \\
 0 & 0 &  \frac{1}{\gamma} z 
  \end{pmatrix}
  \begin{pmatrix}
    \sqrt{\frac{\alpha}{\gamma}} ( \frac{\beta}{ \delta}) I & 0 & 0 \\ 0 &
    \sqrt{\frac{\gamma}{\alpha}} I & 0 \\ 0 & 0 &
    \sqrt{\frac{\gamma}{\alpha}} I  
  \end{pmatrix}
\end{align}
 and by (\ref{TComm}) we can conclude as stated in the following lemma.

\begin{Lemma} \label{Estim}
Let $a$ be in $A,$ $k$ in $A_c$ and $\alpha, \beta, \gamma, \delta
 $ positive real numbers such that $\alpha \beta \leq 1$ and $\gamma
 \delta \leq 1 .$ For $t: = T(a) + k: $ 

\begin{displaymath}
\left \Vert \, \left [ D_{\alpha,\beta} , \pi(t)\right ]\,\right
\Vert  \, \leq \, \max\left \{\frac{\gamma}{\alpha}  , \,
  \frac{\alpha \,\beta^2}{\gamma \, \delta^2}\right \} \left \Vert \, \left [D_{\gamma,\delta} ,
  \pi(t)\right ]\,\right \Vert . 
\end{displaymath}

\end{Lemma}

We will use this result heavily in the computations to come.

\section{The family of compact quantum metric spaces $(A_{\ct},
  L_{\alpha, \beta})$ }

For  a spectral triple $(A, H, D)$ associated to a unital C*-algebra
$\ca,$ Connes has showed that it is
possible to define a metric on the state space $S (\ca)$ of $\ca$ by
the following formula
\begin{equation} \label{distS}
\forall \phi, \psi \in S(\ca):\quad  \dist_{\ca} (\phi, \psi) := \sup
\{\,|(\phi-\psi)(a)| \,|\, \|\, [D,a]\,\| \leq 1\,\}.
\end{equation}
A metric defined in this generality is allowed  to be infinite, but
here we are mostly interested in spectral triples which have the extra
property that the metric defined above is an ordinary metric, which
also is  a metric for the w*-topology on the state space. This aspect
of non commutative geometry has been studied in several articles by
Marc Rieffel \cite{RiLsn} and references there. Rieffel has
generalized this set up to what he calls compact quantum metric spaces.
Here the algebra $A$ of the spectral triple is replaced by an order
unit space and the Dirac operator is not directly present, but
replaced by  a seminorm $L$ on $A. $ In the case where a spectral
triple is present the seminorm is given by $A \ni a \to L(a) := \|\,
[D, a]\,\|.$ Our investigation will not be so
general here, since we will only study degenerations of spectral triples
as constructed in the previous section. On the other hand we will base
our results on Rieffel's memoir \cite{RiMem}, and we will use the
language from that memoir to quantify the impact of the changes of the
parameters $\alpha $ and $\beta.$ We will now recall some definitions
and results  from that memoir. 

\begin{Dfn} 
An order-unit space is a real partially ordered vector space, $A$, with a distinguished element $e$ (the order unit) which satisfies:
\begin{itemize}
\item[(i)] (Order unit property) For each $a\in A$ there is an $r\in \mathbb{R}$ such that $a\leq re$.
\item[(ii)] (Archimedean property) If $a\in A$ and if $a\leq re$ for all $r\in \mathbb{R}$ with $r>0$, then $a\leq 0$.
\end{itemize}
\end{Dfn}
The norm on an order-unit space is given by 
\begin{displaymath}
\Vert a\Vert \,=\, \inf \{ r\in \mathbb{R}: -re\leq a\leq re \}.
\end{displaymath}
Any order-unit space can be realized as a real linear subspace of the
vector space of self-adjoint  bounded operators on a Hilbert space in
such a way that the order unit is the unit operator $I.$ 

\begin{Dfn}
Let $(A,e)$ be an order-unit space, and its dual, $A^{\ast}$. 
The state space  $S(A)$ is defined to be the collection 
of all states, $\mu$, of $A$, i.e. $\mu \in A^\ast$ 
such that  $\mu(e)=1=\Vert \mu \Vert.$ 
\end{Dfn}

\indent Consider now a seminorm $L$ on the order-unit space $(A,e)$
having its  null-space equal to  the scalar multiples of the order unit.
 Then,  for $\mu,\nu \in S(A )$
one can define  a metric, 
$\rho_L$, on $S(A)$ by
\begin{displaymath}
\rho_L(\mu,\nu)\, :=\, \sup\{ \vert \mu(a)-\nu(a)\vert \, | \, L(a)\leq 1\}.
\end{displaymath}
In absence of further assumptions, $\rho_L(\mu,\nu)$  may be infinite.
 It is most often true that the $\rho$-topology on $S(A)$ is finer
 than the weak*-topology. 
\begin{Dfn} \label{Lip}
Let $(A,e)$ be an order-unit space. A Lip-norm on $A$ is 
a seminorm $L$ on $A$ with the following properties:
\begin{itemize}
\item[(i)] For $a\in A$ we have $L(a)=0$ if and only if $a\in \mathbb{R}e.$
\item[(ii)] The topology on $S(A)$ from the metric $\rho_L$ is the weak*-topology.
\end{itemize} 
\end{Dfn}
\begin{Dfn}
A compact quantum metric space is a pair $(A,L)$ consisting of an
order-unit space $A$ with a Lip-norm $L$ defined on it. 
\end{Dfn}

In  our context we have four C*-algebras $C(PH),\, \ca,\, \ct,$ and 
the unitarization of the compacts which we define by
$ \cc := \widetilde{C(PH)} := C(PH) + \bc I_{PH}.$ 
We will now define the order
unit spaces and associated Lip-norms, which we will study.

\begin{Dfn} \label{OUS}
The order unit spaces $A_{\cc}, A_{\ca}, A_{\ct}$ are defined by: 
\begin{align} 
   A_{\cc} \, &:=\ \{ k + \lambda I \, | \, \lambda \in \br, \, k = k^* \in
   A_c \,\} \label{AC}  \text{ Definition \ref{ST} }\\
A_{\ca} \, & := \, \{ a \in A \, |\, a = a^* \, \} \label{AA}\\ 
A_{\ct}\, & := \, T(A_{\ca} ) \, + \,  A_{\cc} \, = \, \{ t \in A_t\,
|\, t = t^*\, \}
 \label{AT}
\end{align}
\end{Dfn}

\begin{Dfn} \label{L}
The seminorms $L_{\cc}, L_{\ca}$ and $L_{\alpha, \beta}$ on $A_{\cc}, A_{\ca}$ and
$A_{\ct}$ are defined  by
\begin{align} \label{LC}
\forall k \in A_{\cc}\cap C(PH)\, \forall  \lambda  \in \br: \, L_{\cc}(k + \lambda I )\, &:= \, \|D_pk\|.\\
\label{LA}
\forall a \in A_{\ca}: \, L_{\ca}(a)\,& := \, \| \,[D,a] \, \| \\ 
 \label{LT}
\forall t \in A_{\ct}:  \, L_{\alpha, \beta}(t)\, & := \, \|
\,[D_{\alpha, \beta},\pi(t)] \, \| .\\
\end{align}
\end{Dfn}
\noindent
The corresponding Minkowski sets or unit balls are defined by

\begin{Dfn} \label{U}
\begin{align}
U_{\cc} \, &:= \, \{x \in A_{\cc} \, |\, L_{\cc}(x) \leq 1 \, \} \label{UC} \\
U_{\ca} \, &:= \, \{ a \in A_{\ca}  \, |\,  L_{\ca}(a) \, \leq \, 1 \,\} \label{UA} \\
U_{\alpha, \beta} \, &:= \, \{ t \in A_{\ct}\, |\,  L_{(\alpha,
  \beta)}( t) \,  \leq \, 1 \,\} \label{UT}
\end{align}
\end{Dfn}
\noindent
The associated metrics are given by 

\begin{Dfn}  \label{dist}
\begin{align}
\forall f, g \in S(\widetilde{C(PH)}): \, \dist_{\cc}(f,g)\, &:=\,
\sup\{|f(k)-g(k)|\,|\, k \in U_{\cc} \,\}  \label{dC} \\
\forall \mu, \nu \in S(\ca): \, \dist_{\ca}(\mu,\nu)\, &:=\,
\sup\{|\mu(a)-\nu(a)|\,|\, a \in U_{\ca} \,\}  \label{dA} \\
\forall \phi, \psi \in S(\ct): \, \dist_{\alpha, \beta}(\phi,
\psi)\, &:=\, \sup\{|\phi(t)- \psi(t)|\,|\, t \in U_{(\alpha,
\beta)} \,\}.  \label{dT}
\end{align}
The diameters of these spaces, which may take the value $\infty $,
 are denoted \newline $\diam_C, \, \diam_{\ca}$  
and $\diam_{\alpha, \beta} $
 respectively.
\end{Dfn}

Let $X$ denote any of the 3 subscripts  $\cc,\,  \ca, \,  (\alpha,
\beta),$ then it follows from \cite{RiStSp} that the metric $\dist_X$
defined above is a metric for the w*-topology on the corresponding
state space,  
if and only if  the
set $U_X$ is separating for the state space, and the image of
$U_X$ in the quotient space $A_X/ (\br I)$ is relatively norm compact
$( A_{\alpha, \beta} := A_{\ct}).$
We start by showing that $L_{\cc}$ is always a Lip-norm and then we
show that any $L_{\alpha, \beta} $ is a Lip-norm if $L_{\ca} $ is so.

\begin{Prop} \label{w*C}
The seminorm $L_{\cc}$ is a Lip-norm.
\end{Prop}

\begin{proof}
Let us first prove that $U_{\cc}$ will separate the state space of
$\widetilde{C(PH)}. $ To this end we remind you that the spectrum
of $D_p$ is discrete and $D_p$ has a compact inverse (on $PH$), 
since it has a
trivial kernel. Let then  $h : = D_p^{-1}$ and it follows  that the set
$\{hyh\, |\, y \in B(PH) \text{ and } y = y^* \} $ is contained in
$A_{\cc} \cap C(PH).$ 
Since $h $ is selfadjoint with trivial kernel ( in $B(PH)$ ) 
this set is norm  
dense in the self-adjoint part of $C(PH),$ and  it 
follows that the set  $\{hyh\, |\, y \in B(PH) \text{ and } y = y^* \}$
 separates
the states and that $U_{\cc}$ will do too. 

Let us define $U_{\cc}^{\circ} \, := \,\{ k \in A_{\cc} \cap C(PH) \,|\, \|D_pk\| \leq
1 \},$ and we will prove that this set is norm compact. From here it
will follow directly that $U_{\cc}/(\br I) = (U_{\cc}^{\circ} + \br I)/(\br I)$
is norm compact.

To see that $U_{\cc}^{\circ} $ is norm compact we will let 
$\eps$ denote a positive real and recall that  $h$ is compact, so 
there exists a finite dimensional spectral projection $E$ for $h
$ such that $\|h(I-E) \| < \eps.$ 
For a $k$ in $U_{\cc}^{\circ} $ we then have 
\begin{align*}
& \|k(I-E)\|\,  = \, \|kD_ph(I-E)\|\, \leq \,  \eps \|kD_p\| \leq \eps.
\\ & \text{ and since  } k \text{ is self-adjoint } \\ 
& \|(I-E)k\|\,  \leq \, \eps .
\end{align*}

\noindent
Then for the set $U_{\cc}^{\circ} $ we get 
$U_{\cc}^{\circ} = U_{\cc}^{\circ}(I-E) + (I-E)U_{\cc}^{\circ}E + EU_{\cc}^{\circ}E,$
 where each operator in
either of the first
two summands is of norm at most $\eps$ and the set $EU_{\cc}^{\circ}E$
 is the unit
ball for some norm on the finite dimensional space $B(EH).$ It then 
follows that $U_{\cc}^{\circ}$ is relatively norm compact. To see that it is  norm
closed, we consider a sequence $(k_n)$
 of elements from $U_{\cc}^{\circ}$ which converges in norm to a compact
self-adjoint operator $k.$ For any spectral projection $E$ of $D_p$ corresponding to
a bounded interval of the real numbers we have that the sequence
$(D_pEk_n)$ is norm convergent with limit $D_pEk,$
and we see that $\|D_pEk\| \leq 1 ,$  and
therefore $k $ belongs to $U_{\cc}^{\circ}.$   

\end{proof}

We will end this section by showing that each of the seminorms
$L_{\alpha, \beta}$ is a Lip-norm
if the seminorm $L_{\ca}$ is a Lip-norm.
This leads to a
detailed study - in the next section - of the two parameter family of
{\em compact quantum metric spaces, } $(A_{\ct}, L_{\alpha,
  \beta}).$ On the other hand we already now
need some estimates on the relations between the various seminorms in
order to prove that $\dist_{\alpha, \beta}$ generates the w*-topology, if
$\dist_{\ca}$ does so.

\begin{Lemma}  \label{Lineq}
Let $\alpha, \beta $ be positive reals such that $\alpha \beta \leq
1,$ $k $ a compact 
 operator in $A_{\cc}$ and
$a$ an operator in $A_{\ca}$ then $t = k + T(a) $ is in $A_{\ct}$ and   
\begin{align}  
 L_{\ca}(a)\, &\leq \, \alpha L_{\alpha, \beta}(t) \label{lin1}\\
 L_{\cc}(k) \, & \leq \, \frac{1 + \alpha \beta }{\beta}L_{\alpha, \beta}(t) \label{lin2}
\end{align}
\end{Lemma}

\begin{proof}
The  first  inequality  follows directly from
(\ref{TComm}) and properties of norms of matrices. 
For the second inequality we use again (\ref{TComm}), the result in
the first inequality and the triangle inequality to obtain 
\begin{align*} 
\beta L_{\cc}(k) \,& \leq \, L_{\alpha, \beta}(t) + \beta \|\,[D_p,
T(a)]\,\| 
\,\leq \, L_{\alpha, \beta}(t) + \beta L_{\ca}(a)\\
& \leq \, (1+ \alpha \beta )L_{\alpha, \beta}(t).
\end{align*} 
and the lemma follows.
\end{proof}

\begin{Prop}  \label{wtop}
If $L_{\ca}$ is a Lip-norm,  then for each pair of positive
reals $(\alpha, \beta)$ such that $ \alpha \beta \leq 
1 $ the seminorm $L_{\alpha, \beta}$ is a Lip-norm.

\end{Prop}
\begin{proof}
To see that $U_{\alpha, \beta}/(\br I)$ is relatively norm
compact we turn back to  Lemma \ref{Lineq}, which 
implies that
for a $t = T(a) + k$ in $U_{\alpha,\beta}$ we have that  
\begin{displaymath} 
a  \in \alpha U_{\ca} \text{ and } k \in
\frac{1+\alpha\beta}{\beta}U_{C},
\end{displaymath} 
so 
\begin{equation} \label{subset}
U_{\alpha,\beta}  \subseteq  \alpha T(U_{\ca}) + 
\frac{1+\alpha \beta}{\beta}U_{C} 
\end{equation}
and 
\begin{displaymath} 
U_{\alpha,\beta}/(\br I)  \subseteq  \alpha T(U_{\ca}/ (\br I)) + 
\frac{1+ \alpha \beta }{\beta}U_{C}/ (\br I). 
\end{displaymath}
Since $\dist_{\ca}$ generates a metric for the w*-topology on $S(\ca)$
the set $U_{\ca}/(\br I) $ is relatively norm compact in $A_{\ca}/(\br I),$ 
and from the proof of Proposition
\ref{w*C} we know that that $U_{\cc}/(\br I) $ is a norm compact
subset of $A_{\cc}/( \br I)$ so
we find that 
$U_{\alpha, \beta}/ (\br I) $ is a relatively norm compact subset of
$A_{\ct}/ ( \br I).$ 
\end{proof}

In the recent article \cite{RiLsn} Rieffel studies Lip-norms which satisfy
some extra conditions, which he needs in order to show certain
results on convergence in the space of compact quantum metric spaces,
equipped with the quantum Gromov-Hausdorff metric. The new seminorms
are called {\em C*-seminorms} and it  seems most likely
that the seminorms we study may possess most of the properties which
a C*-seminorm is required to have. We will not recall all of the
definitions from \cite{RiLsn}, but just recall that one of the
properties is that such a seminorm is demanded to be lower
semicontinuous. In our context this means that the set $\{ t \in
A_{\ct} \, \big| \, \|L_{\alpha,\beta} ( t) \| \leq 1 \} $ is norm
closed. It seems quite unlikely to be the case here since we have
imposed some  regularity conditions in Definition \ref{Tspt}. 
 on the set $A_{\ca}.$  This means that already the seminorm $L_{\ca}$
will probably not in general be lower semicontinuous. On the other
hand we might extend such a seminorm to a larger subalgebra of $\ct$
 and in this way
obtain a lower semicontinuous seminorm, but then it seems difficult
for an operator $t$ in the extended domain for $L_{\alpha, \beta} $ to
control the behavior of the matrix matrix parts of
the commutators of the form $[D_{\alpha, \beta}, t].$ 

We are very thankful to Hanfeng Li, who has showed us, how it is
possible to prove that the seminorm $L_{\alpha, \beta} $ has the two
other properties of a C*-seminorm named {\em spectral stability} and
{\em strongly Leibniz,} provided the original seminorm $L_{\ca} $ has
these properties. On the other hand it seems that the regularity
conditions, we have imposed, may be in conflict with the possibility
for $L_{\ca}$ to be spectrally stable. Mainly inspired by the
classical case we have the impression that the difficluties of this
type may be avoided if we restrict our construction to the special
case, where the domain of definitions for the
seminorms, $A_{\ca} $ and
$A_{\ct} $ are only the  {\em smooth} elements as defined by
Connes in his smoothness axiom of \cite{Co3}. We will present and
discuss this axiom in Section \ref{compacts}.

\section{The compact quantum metric  spaces associated to $\ct$}

In this section we will suppose that the seminorm $L_{\ca}$ is a
Lip-norm, and then by Proposition \ref{wtop} all the tuples $(A_{\ct},
L_{\alpha, \beta } ) $ are compact quantum metric spaces.  This
means that the metric spaces \begin{displaymath}\{\, ( S(\ct), \, \dist_{\alpha, \beta})
\, \big| \, 0 < \alpha \beta \leq 1 \,\} \end{displaymath} 
are equipped with the
w*-topology and hence they are ordinary compact metric spaces. It
seems natural to compare these metric spaces by obtaining Lipschitz
estimates between any pair of two metrics. Based on the Lemma
\ref{Estim} we  can quite easily obtain
such results, which we present just below. The spaces we are studying
are not only compact metric spaces but also {\em compact quantum metric
  spaces } and Rieffel has in the memoir \cite{RiMem}  
developed a distance concept for such spaces called the {\em quantum Gromov-Hausdorff distance.} This
last concept of distance is based on the Hausdorff metric on the
closed subsets of a compact metric space. Gromov has extended this
idea and introduced a distance function defined on pairs of compact
metric spaces, and finally Rieffel \cite{RiMem} has extended Gromov's
ideas to cover the case of compact quantum metric spaces. We will
return to this definition shortly, but first we will treat the
Lipschitz estimates between a pair of metrics $\dist_{\alpha, \beta}
$and $\dist_{(\gamma, \delta )} $ on $S(\ct).$

\begin{Prop} \label{cont}
For any positive reals  $ \alpha,\beta, \gamma, \delta $ such that
$\alpha \beta \leq 1, $ $\gamma \delta \leq 1 $ and any 
 $t$ in $A_{\ct}$
\begin{displaymath}
\min \left \{ \,\frac{\gamma}{\alpha}, \, 
\frac{\alpha\, \beta^2}{\gamma \, \delta^2 \,} \right \} L_{\gamma,
  \delta}(t) \, \leq \, 
 L_{\alpha, \beta}(t) \, \leq \, \max \left \{ \,\frac{\gamma}{\alpha}, \, 
\frac{\alpha\, \beta^2}{\gamma \, \delta^2 \,} \right \} L_{\gamma,
  \delta}(t) 
\end{displaymath}
\end{Prop}
\begin{proof}
We will only prove the right inequality, since the left then follows
by symmetry. 
As usual we have a decomposition of $t$ in $A_{\ct} $  as the sum
$T(a) + k $
 with $a $
in $A_{\ca}$ and $k$ in $A_{\cc} \cap C(PH).$ When going back to the definition in
(\ref{LT}) we get   
\begin{displaymath}
L_{\alpha, \beta}(t) \,= \, \left \| \, \left [ D_{\alpha, \beta}, \pi(t) \right ]\, \right \|
 \end{displaymath}
 and from the results of Lemma \ref{Estim} we 
 then get
\begin{align*}
L_{\alpha, \beta}(t) \,& \leq \,  \max \left \{\,\frac{\gamma}{\alpha}, \, 
\frac{\alpha\, \beta^2}{\gamma \, \delta^2 \,}  \right \}
L_{\gamma, \delta}(t).
\end{align*}
\end{proof}
The results of Proposition \ref{cont} may be applied to the metrics
$\dist_{\alpha, \beta} $ and we can obtain the following proposition. 

\begin{Thm} \label{dcont}
Let $\alpha, \,\beta, \, \gamma, \, \delta $ be positive reals such that $ \alpha \beta \leq 1 $ and $\delta \gamma \leq 1, $ then
the metrics $\dist_{\alpha, \beta}(\cdot, \cdot ) $ and 
$\dist_{\gamma,\delta}(\cdot ,\cdot)
 $ on $S(\ct) $ are Lipschitz equivalent and  satisfy the following
 inequalities 

\begin{align*}
 \forall \phi, \psi \in S(\ct):  \, \, & 
 \min\left \{ \,\frac{\gamma}{\alpha},\, 
\frac{\alpha\, \beta^2}{\gamma \, \delta^2 \,}  \right \} \dist_{\alpha,
  \beta}(\phi,\psi) \, \leq \, \dist_{\gamma,\delta}(\phi,\psi)\\ & \leq \,
\max\left \{ \,\frac{\gamma}{\alpha},\, 
\frac{\alpha\, \beta^2}{\gamma \, \delta^2 \,}  \right \} \dist_{\alpha,
  \beta}(\phi,\psi). 
\end{align*}

\end{Thm}
\begin{proof}
The Proposition \ref{cont} shows that with the notation  from
Definition \ref{U} we get 
\begin{displaymath}
\min \left \{ \,\frac{\gamma}{\alpha},\, 
\frac{\alpha\, \beta^2}{\gamma \, \delta^2 \,}\, \right \} 
 U_{\alpha, \beta} \,
\subseteq \,  U_{\gamma, \delta} \, \subseteq \, 
\max \left \{ \,\frac{\gamma}{\alpha},\, 
\frac{\alpha\, \beta^2}{\gamma \, \delta^2 \,}\, \right \}
 U_{\alpha, \beta}.
\end{displaymath}
The theorem then follows from the definition given in (\ref{dT}).
\end{proof}

We have now seen that any two metrics in this two parameter family of
metrics on $S(\ct)$  are
Lipschitz equivalent, and it follows from this that we can deduce
estimates of the distance with respect to a {\em quantum
  Gromov-Hausdorff metric } between the compact quantum metric spaces
$(A_{\ct}, L_{\alpha, \beta})$ and $(A_{\ct}, L_{\gamma,
  \delta}).$

\medskip
We shall first review, briefly,
 the Gromov-Hausdorff distance for compact metric spaces and Rieffel's
 quantum distance for compact quantum metric spaces.  
We use as references \cite{Li} and \cite{RiMem}.
For any closed subset $Y$ of a metric space $(X, \rho)$ and $r>0$, we denote:
\begin{displaymath}
\mathcal{N}^\rho_r(Y):=\{ x\in X: \exists y\in Y \text{ with }
\rho(x,y)\leq r \}. 
\end{displaymath}  
Let $\mathcal{S}$ denote the class of all non-empty closed subsets of $X$.  The formula, 
\begin{displaymath}
\forall Y,Z \in \mathcal{S}: \text{dist}^\rho_\text{H}(Y,Z):= \inf \{ r: Y\subseteq \mathcal{N}^\rho_r(Z ) \text { and } Z\subseteq \mathcal{N}^\rho_r(Y)\},
\end{displaymath}   
defines a metric (called the Hausdorff metric) on $\mathcal{S}$. One
can also use the notation $\text{dist}^X_\text{H}(Y,Z)$ when there is
no confusion about the metric on $X$. \\ 

\indent Gromov generalized the Hausdorff distance to a distance
between any two compact metric spaces $X,Y$ as follows

\begin{align*}
  &\text{dist}_\text{GH}(X,Y):= \inf \{ \text{dist}^Z_\text{H}(h_X(X),
  h_Y(Y))\big|  h_X:X\rightarrow Z, \, h_Y:Y\rightarrow Z\\ 
  & \text { are isometric embeddings into some compact metric space }
  Z\}.
\end{align*}   

\noindent
One can reduce the space $Z$ above to be the disjoint union $ X
\amalg Y$, and  we shall denote with $\cd (X,Y)$ the set of all
distances 
$\rho $ on  $ X \amalg Y$ fulfilling that the inclusions 
$X,Y \hookrightarrow X \amalg Y$ are isometric embeddings.  
It is then true that
\begin{displaymath}
\text{dist}_\text{GH}(X,Y):=\inf \{  \text{dist}^\rho_\text{H}(X,Y):
\rho \in  \cd (X,Y)\}.
\end{displaymath}

Let $A$ be an order-unit space. By a quotient $(B,\pi)$ of $A$, we mean
an order-unit space $B$ and a surjective linear positive map
$\pi:A\rightarrow B$ preserving the order-unit. Via the dual map
$\pi^\ast :B^\ast \rightarrow A^\ast,$ one may identify $S(B)$
with a closed convex subset of $S(A).$ This gives a bijection between
isomorphism classes of quotients of $A$ and closed convex subsets of
$S(A).$ If $L$ is a Lip-norm on $A$, then the quotient seminorm $L_B$
on $B$, defined by

\begin{displaymath}
L_B(b):=\inf  \{  L(a): \pi (a)=b \}
\end{displaymath}  
is a Lip-norm on $B$, and $\pi^\ast\mid_{S(B)}:S(B)\rightarrow S(A) $ is an isometry for the corresponding metrics $\rho_L$ and $\rho_{L_B}$.

Let $(A,L_A)$  and $(B,L_B)$ be compact quantum metric spaces. The
direct sum $A \oplus B$ has naturally the structure of an order unit space
 with order unit $(e_A,e_B).$ We will let  
$\mathcal{M}(L_A,L_B)$ denote the the set of all Lip-norms $L$ on $A \oplus B$ that induces $L_A$ and $L_B$ under the natural quotient maps 
$ A\oplus B \mapsto A$ and $ A \oplus B \mapsto B $.   
For an element  $L $ in $\cam(L_A, L_B)$  with the associated metric
$\rho_L$   on $S(A \oplus B), $ it is then  possible to consider both of
the compact metric spaces $(S(A), \rho_{L_A})$ and $(S(B), \rho_{L_B})
$ as compact subsets of the compact
metric space $(S(A \oplus B), \rho_L)$, and one can  
 then compute the usual
Hausdorff distance between them. 
This distance is denoted $  \dist_H^{\rho_L}(S(A), S(B)).$ We can
then define a metric on compact quantum metric spaces as follows.
\begin{Dfn} 
Let $(A,L_A)$ and $(B,L_B) $ be compact quantum metric spaces. Then the quantum
Gromov-Hausdorff distance between them is denoted \newline $\dist_q(A,B)$ and
it is defined by \begin{displaymath} 
\dist_q((A,L_A), (B, L_B)) \, := \, \inf \{   \dist_H^{\rho_L}(S(\ca),
S(\cb) ) \, | \, L \in \cam(L_A, L_B)\,\}.
\end{displaymath}
\end{Dfn}
Li gave in \cite{Li} the following description of the Gromov-Hausdorff distance.
\begin{Prop}
Let $(A,L_A)$ and $(B,L_B) $ be compact quantum metric spaces.  Then we have
\begin{align*}
\dist_q((A,L_A), (B, L_B))=& \inf \{ \dist^V_\text{H}(h_A(S(A)), h_B(S(B))): \\ 
& h_A , h_B \text{ are affine isometric embeddings of} \\
& S(A), S(B) \text{ into some real normed space } V \}.
\end{align*}
\end{Prop}

\noindent
This tells us that the quantum Gromov-Hausdorff distance between two
compact quantum metric spaces $(A,L_A)$ and $(B,L_B) $ always will be
larger or equal to the Gromov-Hausdorff distance between the compact
metric spaces $(S(A), \rho_{L_A})$ and $(S(B), \rho_{L_B})$.

Besides Rieffel and Li, there are by 
now several mathematicians who have published articles on
convergence and estimates of distances between compact quantum
metric spaces and even incorporated the extra structure coming from
the theory of operator spaces into their research 
\cite{Ke}, \cite{La}, \cite{Wu} and we have found this very
stimulating for the present work.      

We will now use the results of Proposition \ref{cont} to compute
estimates for the distance between a pair $(A_{\ct},
L_{\alpha, \beta})$ and  $(A_{\ct},
L_{\gamma, \delta})$ of compact quantum metric spaces.
Our construction is based on Rieffel's concept
called a {\em bridge,} but we could not get his concept to fit exactly into
our frame, so we have modified it a bit and incorporated the idea of a
bridge into the proof of the following proposition.  On the other hand our
situation is much simpler than the general situation, considered by
Rieffel,  since the order
unit space is kept fixed as  $A_{\ct}.$

\begin{Prop} \label{dq}
Let $A$ be an order unit space and let $L_1$ and $L_2$ be two
Lip-norms on $A$ for which there exist positive real number $ s < r $ such that 
\begin{displaymath}
\forall a \in A: \quad sL_2(a) \, \leq \, L_1(a)\, \leq \, rL_2(a). 
\end{displaymath}
Define $L_3:= (1/\sqrt{rs})L_1,$ let $\dist_3,\,  \dist_2 $ be the
metrics induced by $L_3, \,L_2$  and let  $\, \diam_3, \,  \diam_2$
 denote the diameters of the compact metric spaces   $(S(A),
\dist_3),$ $ (S(A),
\dist_2)   $ then 
\begin{displaymath} 
 \dist_q((A, L_3), (A, L_2)) \, \leq \,
 \left ( \sqrt{\frac{r}{s}}-1 \right ) \min \{ \, \diam_3, \, \, \diam_2\, \}.
\end{displaymath}

\end{Prop}
\begin{proof}
We first fix an arbitrary {\em base point,  } which in this case means
a state $\sigma$ on $A,$ and then we let $M$ denote an arbitrary
positive real. Later in the argument we will let $M$ increase
unlimited, so you may think of $M$ as a { \em big } positive real. 
We will  let $R$  denote the positive real which is defined by 
\begin{displaymath}
R := \frac{\sqrt{s}}{\sqrt{r} - \sqrt{s}}
\end{displaymath}
and we can then define a seminorm $L$ on $A \oplus A$ by

\begin{align*} 
&\forall a, b \in A: \\
&L(a,b) := \, \max\{\, L_3(a),\, L_2(b), \, RL_3(a-b),\,
RL_2(a-b), \, M|\sigma(a-b)|\,\}. 
\end{align*}
Since  $L$ is defined as a maximum over seminorms, it follows that $L$
is a seminorm on $A \oplus A.$ If $L(a,b) = 0$ then since $L_3$ and
$L_2$ are Lip-norms we see that $a =
\alpha I $ and $b= \beta I$ for some real numbers $\alpha, \beta$ and finally
$\sigma(a-b)= 0$ implies that $\alpha = \beta, $ so $(a,b) =
\alpha(I,I) $ and the first condition for $L$ being a Lip-norm is
established. We will of course also show that $L$ belongs to
$\cam(L_3,L_2),$ and we will address the question of whether $L$ induces
$L_3$ and $L_2$ on the summands first. Let us start by looking at the
first summand and $L_3$ first.  We then define the following sets.
\begin{align*} 
U_L\, &:= \, \{(a,b) \in A \oplus A\, | \, L(a,b) \leq 1 \} \\ 
U_{L|A}\, &:= \, \{a \in A \, | \, \exists b \in A: \, (a,b) \in U_L \} \\
U_2\, &:= \, \{ b \in A\, | \, L_2(b) \leq 1 \} \\
U_3\, &:= \, \{ a \in A\, | \, L_3(a) \leq 1 \}
\end{align*}
In order to prove that $L$ induces $L_3$ it is sufficient to prove
that $U_{L|A} \, = \, U_3,$ so we will do that.
By definition  $L(a,b) \geq L_3(a)$ so for any pair $(a,b) \in U_L $
we have $a \in U_3,$ and then $U_{L|A} \subseteq U_3.$
To establish the opposite inclusion we choose an  $a \in U_3$ and
construct a suitable $b$ such that $(a,b) $ is in $U_L.$ 
  It is a matter of
checking to show that  the element $b$ in $A$ defined by  
 $b\,:= \, \sqrt{s/r}a +(1-\sqrt{s/r})\sigma(a)I$ will do.
  The situation for the second summand is very
 similar, and it turns out that for any $b$ in $A$ such that $L_2(b)
 \leq 1 $ we can define $a$ in $A$ by $a\, :=\, 
\sqrt{s/r}b + (1-\sqrt{s/r})\sigma(b)I,$ and then $L(a,b) \leq 1.$

The seminorm $L$ is defined on all of $A \oplus A $ 
 so  the set $U_L$ will be separating for the states on
$A \oplus A.$ 

We then just have to prove that the set $U_L / (\br (I,I))$
is relatively norm compact in the quotient space $(A \oplus A)/ (\br
(I,I)). $ Let $\tilde\sigma$ denote the state on $A \oplus A$ given by
$\tilde\sigma(a,b) := \sigma(a),$ then it is standard to deduce that 
$U_L / (\br (I,I))$ is relatively norm compact if and only if the set 
$U_{\sigma} \, := \, \{(a,b) \in U \, | \, \tilde\sigma(a,b) = 0 \, \}$
 is  relatively norm compact in $A \oplus A.$ This implies that we may
 define two relatively norm compact sets  in $A$ by  
 $U_{(3, \sigma)} := \{a \in U_3\,|\, \sigma(a) = 0\}$ 
 and $U_{(2, \sigma, M)}:= \{ b \in U_2 \, | \, |\sigma(b) |\leq M \}
. $ For these sets  we find that $U_{\sigma }  \subseteq
 U_{(3, \sigma)}  \oplus U_{(2, \sigma, M)}$ so the metric $\rho_L$ 
 generates the
 w*-topology on $S(A \oplus A).$

We can now use this metric to get an upper  estimate for the quantum
Gromov-Hausdorff distance and we find that for any state $\phi $ on
$A$ 
\begin{align*} 
& \rho_L((\phi,0), (0, \phi))\,  = \,\sup\{\, |\phi(a-b)|\,\big| \, (a,b) \in
 U_L\,\} \\
& \leq \, \sup \{|(\phi - \sigma)(a-b) \big| \, | \, (a,b) \in U_L \, \} +
\frac{1}{M} \\
& \leq \, \min \left \{ \, \frac{\diam_3}{R}, \,  \frac{\diam_2}{R} \, \right \}  + \frac{1}{M}, 
   \text{ since }  (a-b) \in \frac{1}{R}U_3 \cap \frac{1}{R}U_2  
 \end{align*}
By letting $M$ grow we conclude that  
\begin{displaymath} \dist_q((A,L_3), (A, L_2)) \leq \,
  \left (\sqrt{r/s} - 1\right ) \min\{\diam_3, \, \diam_2\, \}.
\end{displaymath}

and the proposition follows.

\end{proof}

\begin{rem} \label{diamrem}
In connection with the proposition  above it may be relevant to note that
the diameters $\diam_2, \diam_3 $ relate in a reciprocal way as the
corresponding seminorms, so we have 
\begin{displaymath}
\sqrt{s/r} \cdot \diam_3 \, \leq \,  \diam_2 \, \leq \,
\sqrt{r/s} \cdot \diam_3. 
\end{displaymath} 
\end{rem}

The special case where the seminorms $L_1$ and  $L_2$ are proportional
is taken out as a corollary.

\begin{Cor} \label{multineq} 
Let $A$ be an order unit space with a Lip-norm $L.$ For any positive
real  $t:$

\begin{displaymath} 
 \dist_q\big( \,(A, L), (A, tL)\,\big) \, \leq \,
   \left  |1-1/t \right |\,\diam_{(A,L)}.
\end{displaymath}
\end{Cor}
\begin{proof}
Suppose $t > 1$ then for the Lip-norm $N := t^2L$ we have 
$L \, \leq \, N \, \leq \, t^2 L. $  The proposition then applies with
$s = 1 $ and $r = t^2, $ so for $tL = (1/\sqrt{sr})N$ we get by the
Remark \ref{diamrem} and the use of the {\em min} option in
Proposition \ref{dq} 
\begin{align*} 
 \dist_q\big( \,(A, L), (A, tL)\,\big) \, \leq &\,
 \left  (\sqrt{r/s} -1 \right )\diam_{tL}\\ =& \,  (t-1)\diam_{tL} \\= & \,(1 -
 t^{-1} ) \diam_L. 
\end{align*} 

For $t < 1 $ and $N = t^2 L$ we get $t^2L \, \leq  \, N \, \leq \, L$
and then 
\begin{displaymath} 
 \dist_q\big( \,(A, L), (A, tL)\,\big) \, \leq \,
 (t^{-1} -1)\diam_L.
\end{displaymath} 

\end{proof}

The corollary above suggests that it could be interesting to see what
will happen for $t$ increasing to infinity, so we will include such a
result. 

\begin{Prop} \label{1point} 
Let $(A, L)$ be an order unit space with a Lip-norm, and let
$(\br, 0)$ be  the one point order unit space with
Lip-norm equal to $0.$ For any positive real $t$ we have the estimate.

\begin{displaymath} 
 \dist_q\big( \,(A, tL),\, (\br,0 ) \,\big) \, \leq \,
    \,\frac{\diam_{(A,L)}}{t}.
\end{displaymath}
\end{Prop}
\begin{proof}
 
We choose and fix a state $\sigma$ on $A$ and let $M$ denote a {\em
  big } positive real. We can then define a seminorm $\hat{L}_t$ on $A \oplus
\br$ by 
\begin{displaymath}
\hat{L}_t (a, s) := \max \{ tL(a),\, M|\sigma(a) - s|\}\end{displaymath}
It is easy to check that $\hat{L}_t$ induces the seminorms $tL$ on $A$
and the zero seminorm on $\br.$ The order unit space  $\br$ has exactly
one state which we denote by $\psi.$ For a state  $\phi$  on $A$ we
can estimate as follows.

\begin{align*}
& \dist_{\hat{L}_t}( (\phi,0), \, (0, \psi)) \, \\ 
&  = \, \sup\{| \phi(a) - s| \, \big| \, \hat{L}_t(a,s) \leq 1 \,\} \\
& \leq \, \sup\{|\phi(a) - \sigma(a)| \, \big| \, tL(a) \leq 1\, \} + \sup\{
  | \sigma(a) - s| \, \big|  \, |\sigma(a) - s |\leq \frac{1}{M}\, \} \\
& \leq \, \frac{\diam_{(A,L)}}{t} + \frac{1}{M}.
\end{align*}
The proposition follows.
\end{proof}

We can then combine some of the results just obtained with Proposition
\ref{cont} to obtain estimates on the variation of the compact quantum
metric spaces $(A_{\ct}, L_{\alpha, \beta}).$ In this connection we
will let $\diam_{\alpha, \beta} $ denote the diameter of this space.  

\begin{Thm}  \label{distThm}
If $\alpha, \beta, \delta, \gamma$ are positive reals such that
$\alpha \beta \leq 1$ and $\gamma \delta \leq 1$ then:

\begin{align*}
&\dist_q\big( (A_{\ct}, L_{\alpha, \beta}),\, (A_{\ct}, L_{\gamma, \delta
})\big) \\ & \leq \,  \left (\,\max \left \{\frac{\alpha\,
    \beta}{\gamma\, \delta},\, \frac{\gamma \, \delta }{\alpha\,
    \beta} \right \} -1 
\, +\,  \left |\, 1-
\frac{\beta}{\delta}\, \right |\, \right )\diam_{\alpha,   \beta},
\end{align*}
\end{Thm}
\begin{proof}
Inspired by Proposition \ref{cont} we define 

\begin{displaymath} 
s \, := \, \min \left \{ \, \frac{\gamma}{\alpha}, \, \frac{\alpha \,
  \beta^2}{\gamma \, \delta^2} \,\right \} \quad r \, := \, \max \left \{\,
\frac{\gamma}{\alpha}, \, \frac{\alpha \, \beta^2}{\gamma \, \delta^2}
\,\right \},
\end{displaymath}
then we get 
\begin{displaymath}
\forall t\in A_{\ct} : \quad s L_{\gamma, \delta}(t) \, \leq \,
L_{\alpha,\beta }(t) \, \leq  r L_{\gamma, \delta}(t).    
\end{displaymath} 
In the notation from  Proposition \ref{dq} 
\begin{displaymath}
\frac{1}{\sqrt{rs}} \, =\, \frac{\delta}{\beta} \, \text{ and } \, 
\sqrt{\frac{r}{s}} \, = \, \max \left \{\frac{\alpha
    \,\beta}{\gamma\,\delta}, \,
\frac{\gamma\,\beta}{\alpha \,\delta}\, \right \},
\end{displaymath} 
  so we have the estimate
\begin{align*} 
&\dist_q \big(\, (A_{\ct}, \, L_{\gamma, \delta} ), (A_{\ct},
(\delta /\beta)L_{\alpha, \beta } )\, \big) \\
&  \leq 
\left (\max \left \{\frac{\alpha\,\beta}{\gamma\, \delta}, \frac{\gamma\,
  \delta}{\alpha \, \beta} \right \} -1 \right )\diam_{\alpha,
  \beta} .
\end{align*} 
We can then use Corollary \ref{multineq} and the triangle
inequality to get  
\begin{align*}
&\dist_q\left (\, \left (A_{\ct},  L_{\gamma, \delta } \right ),  \left ( A_{\ct},
L_{\alpha, \beta } \right )\,\right ) \\& \leq 
\left (\,\max \left \{\frac{\alpha\, \beta}{\gamma\, \delta},\, \frac{\gamma \,
 \delta }{\alpha\, \beta} \right \} -1
\, +\, \left  |\, 1-
\frac{\beta}{\delta}\, \right |\, \right )\diam_{\alpha,   \beta},
\end{align*} 
 and the theorem follows.

\end{proof}

\section{ On limits of $\left (A_{\ct}, L_{\alpha, \beta}\right)$}

In this section  we will keep the set-up from last section so we can  
continue our investigation of the family of
compact quantum metric spaces $\big(A_{\ct}, L_{\alpha, \beta}\big)$
 and study the limiting processes  $\alpha = 1,\, \beta \to 0$ and 
 $\alpha \to 0 , \,
  \beta = 1.$ There are limits in both cases, but they are of
  different nature. In the first case the expression $L_{1,0}$ has a
 an obvious meaning and it follows from (\ref{TComm}) that this will
 be a seminorm on $A_{\ct} .$ This seminorm will be degenerate because
 its kernel will contain all of $C(PH),$  but on the other hand you
 can obtain the seminorm $L_{\ca}$   directly from  $L_{1,0},$ so
 we recover all the ingredients of the original spectral triple via
 this limit process. 
For the family $(\alpha, 1)$ with $\alpha $ decreasing from $1 $ to
$0$ there is no sort of a limit on the level of seminorms, since $\alpha
$ appears in the expression for $L_{\alpha, \beta } $ in the
negative power $1/\alpha,$ but this does not affect the convergence of
the corresponding compact quantum metric spaces since we prove that
the spaces $(A_{\ct}, L_{\alpha, 1 }) $ converge to $(A_{\cc}, L_{\cc}) $ in
the quantum Gromov-Hausdorff metric for $\alpha \to  0.$ 
We have thought of possible interpretations of this result 
and do offer some
remarks concerning the connection to physics in the text below, but we
are not trained physicists, so we are reluctant to make too many comments
 in this direction.

\medskip
The proofs of the results are based on some structural results on the
dual space of a unital C*-algebra.   
Let  the dual space of
$\ct$ be denoted $\ct^*, $ and we will then define two 
subspaces $\cn $ and $\cS$  of $\ct^*$ by 
\begin{align*}
\cn \,  &:= \, \{ \phi \in \ct^* \, |\, \|\phi | C(PH) \|\, = \,
\|\phi\|\,\}\\
\cS\,  &:= \, \{ \phi \in \ct^{\ast} \, |\, \|\phi | C(PH) \|\, = \, 0\,\}
\end{align*}
Here the letters $\cn$ and $\cS$ are chosen because they refer to the
terms {\em normal and singular } functionals on $B(H).$
A priori it is not at all clear that $\cn $ is a subspace, and we will
not prove it here, but recall some results of Effros \cite{Ef} which
are presented just below. For details we refer to Dixmier's book
\cite{Di}  Proposition 2.11.7.

\begin{Prop} \label{l1sum}
With the notation described above, there exist positive
contractive linear projection
operators $N : \ct^* \to \cn $ and $S : \ct^* \to \cS$ such that for
any $\phi$ in $ \ct^*$ 

\begin{align*} 
 N(\phi) + S(\phi) \, = & \, \phi \\
 \|N(\phi)\| + \|S(\phi)\| \, = & \,\| \phi \|
\end{align*}
\end{Prop}

It is easy to identify $\cn$ with the dual space of $C(PH)$ simply by
restricting a functional in $\cn$ to $C(PH).$ The identification the
other way goes via the fact that $B(PH)$ is the second dual of
$C(PH),$
 so the canonical embedding of $C(PH)^*$ into $C(PH)^{***} $
induces an embedding, say $\iota_{\cc},$  of $C(PH)^*$ onto $\cn.$ 

The space $\cS$ may be identified with $\ca^*$ in the following way.
The identification is made  via the homomorphism
$ \rho : \ct \to \ca, $ which was defined in  Definition
\ref{Tspt}. Any  functional $\mu $ in 
$\ca^*$ may be mapped into $\cS$ by the composition $\mu\circ\rho.$
Since the kernel of $\rho$ is $C(PH),$ it follows that this will be an
isometric and order isomorphic mapping of $A^*$ onto $\cS,$ and we
will denote this embedding $\iota_{\ca}.$
 
As an immediate corollary of these identifications 
we get the following result.

\begin{Cor} \label{Pdec}
For any state $\phi$ in $\ct^* $ there exists a unique pair of
  states  $f$ in
  $S(C(PH))$  and $\mu$ in  $S(\ca)$ and a real $\alpha$ in $[0,1] $ such that
$\phi \, = \, (1-\alpha)\iota_C(f) + \alpha \iota_{\ca}(\mu).$  
\end{Cor}

These structures have been studied and generalized in \cite{Al},
\cite{AS} and
in the language of compact convex sets one would say that the two
convex sets \newline
$\iota_{\cc}(S(C(PH)))$ and $\iota_{\ca}(S(\ca))$ form a pair of split
faces of $S(\ct).$

The discussion on how the dual space of $C(PH)$ fits into the dual of
$\ct$ can be applied to the situation when $C(PH)$ is considered as a
subalgebra of $\cc = \widetilde{C(PH)} = C(PH) + \bc I $ too. In this case we
will fix a state $\sigma $ from the space $\cS$ of singular functionals
on $\ct$ and use this state as a basis vector for the one-dimensional 
singular space associated to the decomposition of 
$\widetilde{C(PH)}^* = \cn \oplus \bc \sigma.$ In the general study of
the variation of the metrics on $S(\ct) $ we will use $\sigma $ as a
base point in the w*-compact space  $S(\ct).$

\bigskip 

{\large \bf The limit of $\big(A_{\ct}, L_{1 , \beta}\big)$ as $ \beta  \to 0$}

\bigskip 

This limit is very easy to understand from the point of view of
compact quantum metric spaces. It is simply an affine deformation at
the level of seminorms as it can be seen immediately from the
definitions \ref{ST} and \ref{L}. We will then extend that definition to cover the
pair $(1,0)$ too, and let $L_{(1,0)} $ denote the corresponding
seminorm. We can also still define the unit ball or Minkowski set
$U_{(1,0)} $   for
this seminorm by the definitions given at (\ref{UT}), and it follows
that $A_C \cap C(PH)$ is contained in  $U_{(1,0)}. $ It is then easy to prove
the following result.

\begin{Thm} \label{bto0}
For any $a $ in $ A_{\ca}$ and  $ k $ in $A_{\cc}:$ \begin{displaymath}
  L_{1,\beta}(
T(a) + k)  \, \to \,  L_{1,0}(
T(a) + k) \,  = \,  L_A(a) \, \mathrm{ for } \, \beta \to 0   
\end{displaymath}
\noindent   
For states 
$\phi, \psi$  on $\ct$ with  $\phi = \iota_C(f) + \iota_A(\mu), \,$
$ \psi = \iota_C(g) + \iota_A(\nu)$  the distance formula applied to
the seminorm $L_{1,0}$ gives 
\begin{equation}
\dist_{1,0}(\phi, \psi ) \, = \, \begin{cases}  
0 & \text{ if } \phi = \psi \\
\infty & \text{ if } f \neq g   \\
 \|\mu\| \dist_{\ca}(\mu / \|\mu\|, 
\nu/ \|\mu\|) &\text{ if } f
= g  \text{ and } \mu \neq \nu \end{cases}.
\end{equation}
\end{Thm}

\begin{proof}
Since the kernel of $L_{(1,0)}$ contains all of $A_C,$ it follows from the
distance formula (\ref{distS}) that $\dist_{(1,0)}(\phi, \psi ) = \infty $ if
$f \neq g.$ If $\phi \neq \psi$ and  $f = g$ then $\|f\| = \|g\|  < 1 $ so
$\|\mu\| = \|\nu\| = 1 - \|f\| \neq 0. $ Again the distance formula  and
(\ref{TComm}) give right away that 
\begin{align*}
\dist_{(1,0)}(\phi, \psi) = &\sup\{|(\phi - \psi)(t)|\, \big| \, L_{(1,0)}(t) \leq 1 \,
\}\\ = & \sup\{|(\mu - \nu)(a)|\, \big| \, L_{\ca}(a) \leq 1 \}\\ =& \|\mu\|
\dist_{\ca}(\mu/\|\mu\|, 
\nu /\|\mu\|).
\end{align*}
\end{proof}

We can not prove that the metric distances 
$\dist_{1,\beta}(\phi,
\psi )$ converge to  $\dist_{(1,0)}(\phi,
\psi )$  for $\beta \to 0, $ when the latter is finite, unless we have a trivial extension, but in the cases where the
distance is infinite, i. e. when the 
normal parts, $f$ and $g,$ of the states are different, we
can always give an  estimate of the speed of divergence.

\begin{Prop} \label{dineq}
Let   $0 < \beta \leq 1 $ be a real and  
$\phi, \psi$ states on  $\ct$ with decompositions  $\phi = \iota_C(f) +
\iota_A(\mu), \, 
\, \psi = \iota_C(g) + \iota_A(\nu).$ 
  If $f \neq g$ then there exists a positive real $\gamma$
  such that \newline
 $ \forall \beta \in (0, 1]:\, \dist_{(1, \beta)}(\phi, \psi ) \geq \gamma/ \beta.$
\end{Prop}

\begin{proof}
We will establish a set theoretical inclusion from which the
statement is easy to deduce. 

\begin{equation} \label{Mink1}
\frac{1}{\beta} U_{\cc} \cap C(PH)  \subseteq U_{1, \beta}. 
\end{equation} 
This inclusion  follows from 
 the definitions presented in 
 (\ref{UC}) and the computations which lead to (\ref{TComm}). 
We can then   see that the proposition follows when  we define
$\gamma$ by 
 \begin{displaymath}
 \gamma := \sup\{|(f-g)(k)|\, \big| \, k \in U_{\cc} \cap C(PH)\,\}.
\end{displaymath}
\end{proof}

Suppose $\ca$ is commutative and represents some classical system and
$\ct$ models a quantization of $\ca, $ then for a couple of states on
$\ct$, 
such 
 as $\phi $ and $\psi $  we could look at 
  $f, g$  as their {\em quantum parts } and $\mu, \nu $ as the classical
parts. Then it appears that the limit for $d_{(1, \beta})(\phi,
\psi)$ exists and gives the classical metric, {\em scaled to the size
  of the classical parts }   if and only their quantum
parts are identical.
 Another attempt to make an interpretation is that the inequality
 in the proposition above, implies that in a space
where $\beta $ is small, the quantum parts are far apart; but we do
not want to press this any further right now.

\bigskip

{\large \bf The limit of $\big(A_{\ct}, L_{\alpha, 1}\big)$ as
  $\alpha \to 0$} 

\bigskip

We realized very early on that the family of compact quantum spaces
 $(A_{\ct}, L_{\alpha,1})$ converges {\em pointwise}  as concrete 
metric spaces   towards  \newline  
$(\widetilde{C(PH)}, L_{\cc})$ when $\alpha $ decreases to $0,$ 
 but it took rather long to see that this convergence
actually also works with respect to the quantum Gromov-Hausdorff
metric. Before we prove this result we need a simple estimate.

\begin{Lemma} \label{diamC}
For any  positive functional $f$ in the dual space  $C(PH)^*:$  
\begin{displaymath} 
\sup\{\, |f(k)|\, \big| \, k \in U_{\cc} \cap C(PH)  \, \} \leq \|f\| \diam_C.
\end{displaymath}
\end{Lemma}
\begin{proof}
Let $\eps > 0$ and choose $x $ in  $U_{\cc} \cap C(PH) $ such that $|f(x)|
\geq   \sup\{\, |f(y)|\,|\, y \in U_{\cc} \cap C(PH)  \, \} - \eps/2. $
Since $x$ is compact and $PH$ is of infinite dimension we can find a
positive functional $g$ in $C(PH)^*$ such that $\|g\| = \|f\|$ and
$|g(x)| \leq \eps/2.$ Hence 
\begin{displaymath} 
\|f\|\diam_C \, \geq \, |(f-g)(x)|  \, \geq \,  \sup\{\, |f(y)|\,|\, y \in
U_{\cc} \cap C(PH) \, \} - \eps,
\end{displaymath} 
\end{proof} 

\begin{Thm}   \label{ato0}
For $\alpha, \beta  $ positive reals such that $\alpha \beta \leq 1:$
\begin{displaymath}
\dist_q\left ( \, \left (A_{\ct}, L_{\alpha, \beta} \right ), \, \left
    (A_{\cc}, \beta L_{\cc} \right )\,
\right ) \, \leq \, \alpha\left (\diam_{\ca} + \diam_C \right ).
\end{displaymath}
\end{Thm}

\begin{proof}
We will  define a seminorm on $L$ on $A_{\ct}  \oplus
A_{\cc}$ which induces the given seminorms on each summand.
Let $\sigma$ be a state on $\ct$ which vanishes on $C(PH)$ and let
$M$ be a {\em big} positive real number.  We
 can then define the  seminorm $L$.  
\begin{align*}
&\forall a \in A_{\ca} \, \forall k, h \in A_{\cc} \cap C(PH)
 \, \forall s \in \br:
\,\,  L((T(a)+k, h + s I))\,   :=\, \\ &\max  \{\, L_{(\alpha,\beta)}
(T(a)+k),\,
\beta L_{\cc}(h),  
\frac{1}{\alpha} L_{\ca}(a),\, \frac{1}{\alpha}L_{\cc}(k-h),\,
M|\sigma(T(a) - sI)|\,\}  
\end{align*}

Let us show that the seminorm induced by $L$ on $A_{\ct}$ is $L_{(\alpha,
\beta)}.$ By definition we always have $L((T(a)+k, h + s I)) \geq
L_{(\alpha,\beta)}(T(a)+k)$ so it is enough  to prove that for a given $t \, =
\, T(a) + k$ with $a$ in $A_{\ca}$ and $k $ in $A_{\cc} \cap C(PH)$ we can
find  an  $h $ in
$A_{\cc} \cap C(PH)$ and  an $s $ in $\br$ such that $L((T(a)+k, h + s I)) =
L_{(\alpha,\beta)}(T(a)+k).$ We will prove that   $h :=
(1 + \alpha \beta)^{-1}k$ and $s:= \sigma(T(a))$ will work. To this end we
may  without loss of
generality assume that $L_{\alpha, \beta}(T(a)+k) =1 , $ 
and then by    (\ref{lin1})
it follows that $L_{\ca}(a) \leq \alpha, $ and by (\ref{lin2}) we find
that  $L_{\cc}(k) \leq (1+\alpha \beta)/\beta.$ From here it is
easy to prove that $L(T(a)+k, h + sI ) = 1.$  
For the seminorm induced by $L$ on $A_{\cc}$ we also get by
definition that  $L((T(a)+k, h + s I)) \geq
 \beta L_{C}(h).$ Let then an $h + s I $ be given in   in $A_{\cc} $ and
define $a := s I, $  $k:=h,$ then it is again a matter of computation
to show that $L((T(a)+k, h + s I)) = \beta L_{C}(h).$ \\
\indent We will then show  $\dist_q\big(  (A_{\ct}, L_{(\alpha, \beta )}), \,
(A_{\cc}, L_{\cc})  \leq  \alpha (\diam_{\ca} + \diam_C) $
 by showing that for each positive $\eps$ and any 
 state $\phi $ on $\ct$ there exists a state
 $\psi$ on $\widetilde{C(PH)}$ such that for the metric $\rho_L$ on the
 state space  of $A_{\ct} \oplus A_{\cc}$ we have $\rho_L((\phi,0),
 (0,\psi)) \leq  \alpha (\diam_{\ca} + \diam_C) + \eps , $ and vice versa.\\
 \indent For a state $\phi $ on $\ct$ we can write $\phi = \iota_{\cc}(f) +
\iota_{\ca}(\mu) $ for positive functionals $f$ on $C(PH)$ and $\mu $
on $\ca.$ Let $\hat f$ denote the extension - with the same norm - 
of $f$ to
$\widetilde{C(PH)}$,
then  the functional $\psi$ is defined as $\hat f + \| \mu \| \sigma$ on
$\widetilde{C(PH)}$ and we get the following string of inequalities 
\begin{align*}
&\rho_L(( \phi,0),\, (0, \psi))\,\\  = \,& \sup\{ \, |\phi(T(a)+k) - \psi(h
+sI)\, \big| \,\, |\, L((T(a)+k, h + s I))\, \leq \, 1\,\} \\
 \leq \,& \sup\{\,|\phi(T(a)) - \sigma(T(a))| \, \big| \, L_{\ca}(a) \leq
\alpha\}  \\ + &\sup\{\,|\sigma(T(a)) - s |\, \big| \, |\sigma(T(a)) - s|
\leq 1/M\} \\ +  &\sup\{\,|f(k-h)|\, \big| \, L_{\cc}(k-h) \leq
\alpha\} \text{ which  by Lemma \ref{diamC} }\\ 
 \leq \, & \alpha ( \diam_{\ca} + \diam_C) + \frac{1}{M},
 \end{align*}
Given  a state $\psi $ on $\widetilde{C(PH)}$ we can write $\psi = \hat f +
(1-\|f\|) \sigma $ for a positive functional $f$ on $C(PH)$ of norm at
most 1. Then the functional $\phi$ is defined as 
$\iota_C(f) + (1-\| f \|) \sigma $ on $\ct,$ and we get as above.  
\begin{align*}
&\rho_L(( \phi,0),\, (0, \psi))\,\\  = \,& \sup\{ \, |\phi(T(a)+k) - \psi(h
+s I)\,| \,\, \big| \, L((T(a)+k, h + s I))\, \leq \, 1\,\} \\
 \leq \,& \sup\{\,|\phi(T(a)) - \sigma(T(a))\,|  \, \big| \, L_{\ca}(a) \leq
\alpha\}  \\ + &\sup\{\,|\sigma(T(a)) - s |\, \big| \, |\sigma(T(a)) - s|
\leq 1/M \} \\ +  &\sup\{\,|f(k-h)|\, \big| \, L_{\cc}(k-h) \leq
\alpha\} \\ 
 \leq \, &\alpha ( \diam_{\ca} + \diam_C) + \frac{1}{M},
 \end{align*}
 and the theorem follows.  
\end{proof}
The inequalities just above show, that when $\alpha \to 0$ then
the system seems to forget how it was created and only the very basic
structure of the  quantum infinitesimals modelled by $C(PH)$ are left
visible.

\section{ A quantum metric on the set of parameters  $\cp \, : = \, \{(\alpha, \beta) \in
  \br^2\,| \,  \alpha \geq 0,\, \beta > 0, \,  \alpha \beta \leq 1\,
  \} \cup \{(0,\infty)\}$}   

The quantum Gromov-Hausdorff metric on our two-parameter family of
compact quantum Hausdorff spaces naturally define a metric on the
parameter space, say $\cp^{\circ} := \{(\alpha, \beta) \in
  \br_+^2\,| \,  \alpha \beta \leq 1\, \}, $ 
and we want to get an impression on the sort of
metric space we can obtain this way. We have not made a very detailed
study of this but we show that some balls in this metric are unbounded with
respect to the Euclidian distance in $\br^2.$ We also show, the other
way around, that some sets which are bounded with respect the Euclidian
metric are unbounded with respect to the {\em quantum-}metric.  
Based on the results in Theorem
\ref{ato0} we realized  that it is reasonable to extend the
parameter space to the space $\cp,$ defined below. 

\begin{align*}
\cp\,& := \, \{(\alpha, \beta) \, |\, \alpha \geq 0 , \, \beta > 0 , \,
\alpha \beta \leq 1 \,\} \cup \{(0, \infty)\} \\
\text{ for }0 \, < \,\beta < \infty &: \  (A_{0,\beta},\,
L_{0,\beta})\,: =\, (A_{\cc},\, \beta L_{\cc}) \\
 (A_{0,\infty}, L_{0,\infty}) \,& :=\, (\br, \, 0).
\end{align*}

The points we have added are also compact quantum metric spaces, and
it turns out that they fit in very well with respect to the quantum
Gromov-Hausdorff metric.

\begin{Prop}  \label{CPset}
Let $\beta_0 > 0$ then the subset $\cp_{\beta_0} \, := \, \{(\alpha,
\beta ) \in \cp \, |\, \beta \geq \beta_0 \,\}$  is  compact with
respect to the metric inherited from the quantum Gromov-Hausdorff
distance. 
\end{Prop}
\begin{proof} 
Fix a positive  $\eps$  and define $\beta_1 \, := \, \max \{\beta_0,
2(\diam_A + \diam_C)/ \eps\}.$ For any pair $(\alpha, \beta) $ in $\cp$
with $\beta \geq \beta_1$ we get $\alpha  \leq \beta_1^{-1}$ and by
Theorem \ref{ato0} 
\begin{displaymath} \dist_q((\alpha, \beta ), (0,\beta) ) \, \leq \,
  \alpha( \diam_A + \diam_C) \,\leq \, (\diam_A + \diam_C )/\beta_1.
\end{displaymath}
By Proposition \ref{1point} 
\begin{displaymath} 
\dist_q((0,\beta), \, ( 0, \infty) ) \, \leq \, \diam_C/ \beta \, \leq
\, \diam_C/ \beta_1.
\end{displaymath} 
hence it follows that for $(\alpha, \beta ) $ in $\cp$ with $\beta
\geq \beta_1$ this point is in the ball of radius $\eps$ with centre
in $(0, \infty ).$ We are then left with the set $\{(\alpha, \beta)
\in \cp_{\beta_0} \, | \, \beta_0 \leq \beta \leq \beta_1\,\}$ and we
will divide this set into two sets dependent on a positive real
$\delta $ which we define by 
\begin{displaymath} 
\delta \, := \, \min \{ \frac{\eps}{3(1+ \diam_{\ca} + \diam_C)},
\frac{1}{\beta_1} \} 
\end{displaymath}
and the sets become 
\begin{align*}
\cx \, & := \, \{(\alpha, \beta)
 \, | \, 0 \leq \alpha \leq \delta \text { and }
\beta_0 \leq \beta \leq \beta_1  \, \}
\\
\cy \, & := \, \{(\alpha, \beta)
\, | \, \delta  \leq \alpha  \text { and } \beta_0 \leq
\beta \leq \beta_1 \text{ and } \alpha \beta \leq 1\, \}
\end{align*}

\smallskip
\noindent
 The results from Theorem \ref{distThm} 
show that the usual Euclidian metric and  the metric $\dist_q$
generate the same topology on the subset
$\cy,$ so this set is compact. For the set $\cx$ we can look at the
subset $\cz$ which we define by 
\begin{displaymath} 
\cz := \{(\delta, \beta) \, |\, \beta_0 \leq \beta \leq \beta_1 \,\}
\end{displaymath}
Since $\cz$ is also a subset of $\cy,$ it is compact for the quantum
metric $\dist_q,$
 by the result above, and we can find a finite number of points
$\{( \delta , \beta_i)\, | \, i \in J \} $
 in $\cz$ such that any point in $\cz$ is
within distance $\delta$ from a point of the form $( \delta ,
\beta_i).$ For any point   $(\alpha, \beta )$ in  $\cx, $ we get from
Theorem \ref{ato0} that $\dist_q((\alpha, \beta), (0, \beta)) \leq
\eps/3$ and for suitable $\beta_i$ we get 
\begin{align*} 
 \dist_q((\alpha, \beta), (\delta, \beta_i)) \, \leq \,& \dist_q((\alpha,
 \beta), (0,  \beta))\\  + &  \dist_q((0, \beta), (\delta, \beta)) +
 \dist_q((\delta,  \beta), (\delta, \beta_i)) \leq \eps,  
 \end{align*}
 and the
proposition follows.   

\end{proof}

We will then look at the subsets of $\cp$ such that $\alpha \geq
\alpha_0.$ Here the situation is quite the opposite since these sets
will be unbounded with respect to the quantum metric on $\cp. $  To
see this we fix a positive $\gamma \leq 1$ and we will study behavior
of the metric along the
hyperbola $\ch_{\gamma} \, := \{(\alpha, \beta) \in \br^2_+\, |\, \alpha
\beta = \gamma \,\}.$ We see that   
 the seminorms corresponding to the points on $\ch_{\gamma} $
  are all proportional and for any positive real $s$ we see   from the
  Definition \ref{L}  $L_{(\gamma / s), s} = s L_{\gamma, 1}, $ so
  the space is well understood along each of these curves. In
  particular, for  the diameters we have $\diam_{(\gamma/s), s} =
  (\diam_{\gamma, 1 }) /s,$ so for $s \leq 1 $ and $s$ decreasing to
  $0,$ we get immediately the following estimate.
  
\begin{Prop} 
For positive reals $\gamma, s$ such that  $0 <\gamma \leq 1$ and $ 0
< s \leq 1/\gamma :  $
\begin{displaymath} 
\dist_q((A_{\ct}, L_{\gamma, 1}), (A_{\ct}, L_{(\gamma/s),
  s}))\  \geq (1/2)(s^{-1}-1) \diam_{\gamma, 1} 
\end{displaymath} 
\end{Prop}

\begin{proof}
Let $\eps > 0$ and let $\delta := \dist_q((A_{\ct}, L_{\gamma, 1}),
(A_{\ct}, L_{(\gamma/s), s})).$ For a pair of states, say  $\phi,
\psi  $ on $\ca_{\ct} $ such that $d_{(\gamma/s),s}(\phi, \psi) \geq
(\diam_{\gamma, 1})/s - \eps/3 $ we can find approximating states -
with respect to 
 $(A_{\ct}, L_{\gamma, 1}),$  - say
$\mu $ and $\nu$  on
$\ca_{\ct}$ such that

\begin{align*}
(\diam_{\gamma, 1 })/s - \eps/3  \leq & \dist_{(\gamma/s),  s}(\phi,
\psi) \\
  \leq & 2\delta +  2 \eps/3 + \dist_{\gamma,1} (\mu, \nu) \\ \leq &
  2
  \delta + 2 \eps/3 +
  \diam_{\gamma, 1},
\end{align*}
and the proposition follows.
\end{proof}

For the vertical intervals $\{ (\alpha, \beta) \, \big| \, 0 < \beta
\leq \alpha \,\}$ we get that they are all unbounded with respect to
this new metric. This follows easily from Proposition \ref{dineq}, and
we will state it formally in the following proposition. 

\begin{Prop} \label{UBset}
For a fixed $\alpha_0 > 0$ there exists a positive $\gamma $ such that
 \begin{displaymath} \forall \beta \in \, ]0, \, 1 /  \alpha_0]: \,
   \diam_{\alpha_0, \beta} \geq \gamma/ \beta . \end{displaymath}
\end{Prop} 

\section{Applications to the compacts and an investigation of Connes'
  7 axioms for this spectral tripæle.} \label{compacts}

Right after the Definition \ref{Ttype} of the Toeplitz extension of a
C*-algebra $\ca,$ we remarked that it is
debatable if the generalized Toeplitz algebra should be defined as the
C*-algebra generated 
by $P\ca|PH$ alone or - {\em as we have chosen } -  
the one generated by this set plus the
compacts. The difference is a trivial extension, but for the algebra
$\widetilde{C(H)}$  it is a rather crucial difference, when this algebra is
considered to be a trivial extension of the one dimensional C*-algebra
$\bc I.$  Our first example here shows
that our construction offers a variety of spectral triples for the
unitarized 
compacts. On the other hand, for the Podle{\`s} sphere our construction
gives an algebra which has more compacts than the universal C*-algebra
for the Podle{\`s} sphere has. If we just had used the C*-algebra
generated by $P\ca|PH,$ we would have obtained the right algebra here.

\begin{exmp} 
Let $H$ be a separable infinite dimensional Hilbert space and let $\ca
:= \bc I $ be the unital C*-algebra generated the unit $I$ on $H$. Let
$D$ be an unbounded  self-adjoint invertible operator on $H$ with
compact inverse and let the projection $P := I.$ We now have a
spectral triple $(\ca, H, D)$ and a quadruple $\big((\ca, H, D), P \big)$ of
Toeplitz type. Our construction will then give a C*-algebra $\ct : =
\widetilde{C(H)},$ a Hilbert space $K:= H \oplus H$ and a representation
$\pi$ of  $\ct$ on $K$  by 

\begin{displaymath} 
\forall k + \gamma I \, \in \, C(H) + \bc I \quad \pi(k + \gamma I
)\,:=\, \begin{pmatrix} k + \gamma I & 0 \\ 0 & \gamma I \end{pmatrix}.
\end{displaymath} 

\noindent
The Dirac operator then becomes 

\begin{displaymath} 
D_{\alpha, \beta} \,:=\, \begin{pmatrix} 0 & \beta D \\ \beta D &
  \frac{1}{\alpha} D \end{pmatrix}.
\end{displaymath} 
\end{exmp}

You may notice that the part  $1/ \alpha D$ has
no effect, for this spectral triple and this leads to the following
proposition which will yield many more spectral triples associated to
$\widetilde{C(H)}.$

\begin{Prop} \label{CTrip}
Based on the notation in the example above let $T$ be an unbounded
densely defined and closed operator on $H.$ If  $|T| $ is invertible
with compact resolvent then for 
\begin{displaymath} D  \,:=\, \begin{pmatrix} 0 & T^* \\ T &
  0 \end{pmatrix},
\end{displaymath}
the set $\big(|T|^{-1}C(H)|T|^{-1} + \bc I , K , D \big)$ is a spectral triple
associated to $\widetilde{C(H)}.$
\end{Prop}

\bigskip

In the article \cite{Co3} Connes lists 7 axioms for Non Commutative
Geometry. In the present case of a spectral triple associated to the
unitarized compacts $\widetilde{C(H)} $  the dimension must be $0$ so
the dimension is even and for the Dirac operator $D$  the growth of the
eigenvalues must be such that for any positive real $s$ the operator
$(I + D^2)^{-s/2} $ is of trace class. There should also be a grading
$\gamma$
and a conjugate linear operator $J$ which relate in certain ways. We
can provide candidates for these ingredients, which seem natural to
us, but they will not fulfill all of Connes' axioms. We will therefore
present the candidates for $D, \, \gamma, \, J,  $ check each of the
axioms and show what sort of problems we are facing. We keep the
notation from above in this section and define 

\begin{Dfn} 
\begin{itemize}
\item[]

\item[(D)]

Let $T$ be a self-adjoint unbounded operator  with trivial kernel,
such that for any
positive real $s$ the operator  $|T|^{-s}$ is of trace class,
then the Dirac operator $D$ is defined on $H\oplus H $ by
\begin{displaymath}
 D\, : = \,
\begin{pmatrix} 0 & T \\ T & 0 \end{pmatrix}. \end{displaymath}

\item[($\gamma$)] The obvious choice for this unitary seems to be the
  unitary on $H \oplus H, $ given by  
\begin{displaymath}
 \gamma \, : = \,
\begin{pmatrix} I & 0 \\ 0 & -I \end{pmatrix}. \end{displaymath}

\item[(J)] It is not so obvious what to choose here,
  since our setup is not the same
  as the one Connes clearly has in mind. In \cite{Co3} Connes obtains the $J$
  operation from a standard representation of a self-adjoint algebra
  of bounded operators.  This is not what we have here for $C(H),$ but
  we have anyway a candidate for $J$ which seems reasonable. First we
  define $j : H \to H$ by chosing an orthonormal basis $(\xi_n)$ for
  $H$ consisting of eigenvectors for $T.$ Then we define $j$ on $
  \lambda \xi_n$ as $ \bar{\lambda} \xi_n$ and extend this to a
  conjugate linear isometry of $H$ onto $H.$ The choice for $J$ is
  then given by
\begin{displaymath}
 J \, : = \,
\begin{pmatrix} 0 & j \\ j & 0 \end{pmatrix}. \end{displaymath}
Remark that 
\begin{displaymath}
 J \pi(a^* + \bar\lambda I) J \, = \, \begin{pmatrix} ja^*j  +   \lambda I & 0
   \\ 0 & \lambda I \end{pmatrix}, \end{displaymath}
so for any operators $a + \lambda I $ and $b + \mu I $ the operators
$\pi( a + \lambda I)$ and $J\pi(b^* +\bar \mu I)J$ do commute. 
\end{itemize}\end{Dfn}

\noindent
We will then look at the 7 axioms taken from \cite{Co3}  one by one,
but first we will define the algebra $\ca$ of smooth elements as the
operators $ a \in \widetilde{C(H)} $  such that for $\delta(x) : =
[|D|,x ] $ we have for any $a \in \ca$ and any natural number $m$ 
both $\pi(a) $ and $[D, \pi(a)] $ are in the domain
of $\delta^{m}.$  Inside this algebra $\ca$ we have a norm dense
subalgebra of operators of finite rank which we denote $\ca_0.$ This
algebra is defined via the orthonormal basis  $(\xi_n),$ from above,
consisting of eigenvectors for $T.$ The algebra $\ca_0$ is then the
linear span of the matrix units $a_{ij} := \langle ., \xi_j\rangle
\xi_i.$  
\begin{itemize}
\item[(1)] {\em The operator $D^{-1}$ is an infinitesimal of infinite
    order. } This is fulfilled by the demand that $|T|^{-s} $ is of
  trace class for any psoitive real $s.$
\item[(2)] For any pair of elements $a, b \in \ca: [[D,
  \pi(a)],J \pi(b^*)J] = 0 .$  This demand can not be met, but we have -
  as in Dabrowski's paper \cite{Da} -  $ \forall a, b \in \ca: [[D,
  \pi(a)], J\pi(b^*)J] \in  C(H \oplus H). $ 

\item[(3)]  {\em Smoothness } The smoothness axiom is fulfilled by the
  definition of the algebra $ \ca.$ 

\item[(4, 5, 6)] We can not show that the spectral triple we
  investigate fulfills any of these 3 axioms. We only have the grading
  $\gamma$ and it seems to be uniquely determined by its basic
  properties. 
\item[(7)] We have an operator $J$ such that $[\pi(a), J\pi(b^*)J] =
  0,$  but $\gamma$ and $J$ do not fit with the  {\em reality properties } of the table in \cite{Co3}
  p. 162. We get  $J^2 = I, \,
  JD = DJ, \, J\gamma = - \gamma J.$ For $n=0$ the first two
  identities are as in the table, but the last one should have been
  $J\gamma = \gamma J.$   
 
\end{itemize}

We will now turn to the Podle{\`s} sphere, $C(S^2_{qc}),$
  \cite{Po} and base the
presentation here on Section 4 of that paper and on  Chakraborty's
description, \cite{Ch},  of a concrete faithful representation for
this algebra. Chakraborty's purpose was to some extent the same as
ours since he wanted to create a  Lip-norm for the Podle{\`s} sphere,
based on the fact that this C*-algebra is an extension of the
classical Toeplitz algebra by the compacts.  In the latter
presentation the parameter $\mu$ from Podle{\`s} paper is replaced by
the letter $q,$ which now seems to be standard, so we will use this
notation.   

\smallskip

For $c, q $ reals such that $c > 0$ and $0 < |q|
< 1 $  the Podle{\`s} sphere, $C(S^2_{qc})$ 
 is the universal C*-algbera generated by 2 operators $A$ and $B$
 which satisfy the following relations:

\begin{displaymath} 
A = A^*, \, \, BA = q^2AB, \, \, B^*B = A - A^2 + cI, \, \, BB^* =
q^2A  - q^4I + cI.
\end{displaymath}

Let $\ct(\bt)$ denote the classical Toeplitz algebra for the unit
circle and let $\rho : \ct(\bt) \to C(\bt)$ denote the canonical
surjective homomorphism. All the algebras  $C(S^2_{qc})$  turns out to be
isomorphic \cite{Sh} and  can  be
described by 
\begin{displaymath} 
C(S^2_{qc}) \, = \, \{(x, y ) \in  \ct(\bt) \oplus  \ct(\bt)\, |\,
\rho(x) = \rho(y)\}.
\end{displaymath}   

By 
So we can see that $C(S^2_{qc}) $ is an extension of the Toeplitz
algebra by the compacts. 

Let us consider the standard spectral triple associated to the
C*-algebra $\ca$ of 
continuous functions on the unit circle. For the algebra $A$ we can take
the functions whose Fourier coefficients form  rapidly
decreasing sequences and the Dirac operator is 
$\frac{1}{i}\frac{d}{d \theta}.$  
As the Hilbert space $H$ we take $L^2(\bt)$ and the
projection $P$ is the projection onto $H_+,$ the closed linear span of
the eigen functions 
$ e^{i n \theta}, \, n > 0,$ corresponding to the positive eigen
values. This is not the usual definition, where the constant function
$I$ usually is assumed to be in $H_+.$ This will not change the
construction {\em qualitatively } but it will have the nice consequence,
that the restriction of $D$ to $H_+$ is invertible with a compact
inverse. Finally we let $H_-$ denote the orthogonal complement of
$H_+$ in $H.$ 
  The quadruple $\big((\ca, H, D ), P \big) $ is then 
  of Toeplitz type and the Definition \ref{ST} gives a spectral
  triple $(A_t, K, D_{\alpha, \beta})$ for the ordinary Toeplitz
  algebra. Recall that the Hilbert space $K$ is given as $H_+ \oplus
  H_+ \oplus H_-$, so we can define a projection $Q$ of $K$ onto the
  first two summands  and it follows from Definition \ref{ST},(v),
  that $D_{\alpha, \beta}$ commutes with $Q$. By checking the same
  definition's point (iv) it can be seen that for any $t$ in $\ct(\bt)$
  the commutator $[\pi(t), Q]$ is compact since it is nothing but
  the embedding of the operator $[\rho(t), P]$ into
  $B(K).$  In order to show, that we now have a quadruple of Toeplitz
  type in the set
  $\big((A_t, K, D_{\alpha, \beta}), Q \big),$ we then, according to
  Definition \ref{Tspt} only have to prove that $D_{(\alpha,
    \beta)Q}:= D_{\alpha, \beta}| QK$ has trivial kernel, but this
  follows easily from the description of $D_{\alpha, \beta} $ given
  in Definition \ref{ST} point (v) and the fact that $D_P$ is assumed
  to be injective.    As mentioned
  above we will not consider the extended algebra as defined in
  Definition \ref{Ttype}. Instead we will define $\ct_{\ct}$ as the
  C*-algebra on $H_+ \oplus H_+$ generated by $Q\pi(t)|QK.$ It is not
  difficult to see that this C*-algebra is exactly the one which
  above is described as the algebra $C(S^2_{qc}).$ Our construction
  can give a family of spectral triples associated to $C(S^2_{qc}) +
  C(H_+ \oplus H_+), $ and let $ (A_{tt}, K_{tt}, D_{tt}) $ denote
  such a set. Then  we are left with the problem to realize how the
  algebra $A_{tt}$ relates to the sum $C(S^2_{qc}) + C(H_+ \oplus H_+).$
  To deal with this question
 we first remark that by Definition \ref{ST}  point (ii) any
  element in  $A_{tt}$ is a sum of an element related to a {\em
    differentiable symbol }  and {\em a differentiable compact.}
  Consequently we only have to see how a {\em differentiable
    compact, } say $C$ in  $C(H_+ \oplus H_+))$  
behaves with respect to the splitting as a
  sum of a diagonal operator and an off diagonal operator,

 \begin{displaymath} 
C \, = \,\begin{pmatrix} u & v \\ x & y \end{pmatrix}\, =\,
\begin{pmatrix} u & 0\\ 
  0 & y \end{pmatrix}\, +\, \begin{pmatrix} 0 & v \\ x & 0 \end{pmatrix}.
\end{displaymath}

By Definition \ref{ST} point (i) the matrix above is a {\em differentiable
compact} if and only if both of the products $D_{(\alpha, \beta)Q}C$
and $CD_{(\alpha, \beta)Q}$ are bounded and densely defined. In matrix
forms these products are as seen below
\begin{align*} 
D_{(\alpha, \beta)Q}C \, &= \, \begin{pmatrix}\beta D_P x & \beta D_P
  y \\ \beta D_P u + \frac{1}{\alpha}D_P x & \beta D_P v +
  \frac{1}{\alpha}D_P y \end{pmatrix} \\ 
CD_{(\alpha, \beta)Q} \, &= \, \begin{pmatrix}\beta v D_P  & \beta u
  D_P + \frac{1}{\alpha} v D_P \\ 
\beta y D_P & \beta x D_P + \frac{1}{\alpha} y D_P  \end{pmatrix}.  
\end{align*}
From here it follows that the products $D_{(\alpha, \beta)Q}C$
and $CD_{(\alpha, \beta)Q}$ are bounded and densely defined if and
only if all of the operators $u, v, x, y $ belong to the algebra $A_c,$
as defined in Definition \ref{ST} point (i). This has the immediate
consequence that if we replace $A_{tt} $ by $A_{ttd}$ which we define
by 
\begin{displaymath}
A_{ttd} \, := \, \left \{ \begin{pmatrix} a & 0 \\ 0 & d \end{pmatrix}
  \, \big  |
\,    \begin{pmatrix} a & 0 \\ 0 & d \end{pmatrix} \, \in \, A_{tt} \,
 \right \}, 
\end{displaymath}
then $(A_{ttd}, K_{tt}, D_{tt})$ is a spectral triple associated to the
universal C*-algebra for the Podle{\`s} sphere.

\begin{rem}
We have been asked by the referee, if there are some connections
between our example of a spectral triple for the Podle{\`s} sphere 
and the ones obtained in \cite{DDLW} and \cite{DS}.
 We have, but  in vain,
 tried to answer this
question, and it is our impression that there is no simple connection.  
\end{rem}

\section{Odd and even extensions an analytic K-homology}

Our constructions in this paper produce odd spectral triples and this
seems not to be the right setup for algebras containing the compact
operators. It is quite easy to produce an even spectral triple from an
odd one by doubling the representation and introduce the Dirac
operator $\hat{D} $ on the Hilbert space $K \oplus K$ which is given
by \begin{displaymath} \hat{D} = \begin{pmatrix} 0 & D \\ D & 0
  \end{pmatrix}
\end{displaymath}
The grading $\gamma$ is then given on $K \oplus K$ by $\gamma(\xi,
\eta ) := (\xi, - \eta).$ 
We could have performed all our computations in this setting, but it
would not give any new insights with respect to the metric properties we
have been investigating in this article, so we have
not pursued a presentation this way.  

Extensions of unital C*-algebras by the compacts as we do it here is
described in Higson and Roe's book \cite{HR} Chapter 5. So according
to that description an extension of the sort we are looking at
corresponds to an element in the reduced analytic K-homology of the
unital C*-algebra $\ca.$ But our construction is only designed for
projections $P$ in the commutant of $D.$


\begin{thebibliography}{9999}


\bibitem[1]{Al} E. M. Alfsen. {\em Compact convex sets and boundary
    integrals.} Springer-Verlag, 1971. 

\bibitem[2]{AS} E. M. Alfsen, F. W.  Shultz. {\em State spaces of operator algebras. Basic theory, orientations and C*-products}, Birkhäuser, 2001.

\bibitem[3]{BJ}
S. Baaj, P. Julg,  {\em Th{\'e}orie bivariante de Kasparov et
  op{\'e}rateurs non born{\'e}s dans les $C^*$-modules hilbertiens}, 
C.R. Acad. Sci. Paris, Serie I, {\bf 296} (1983), 875--878. 


\bibitem[4]{Be} I. Belgradek, {\em Degenerations of Riemannian manifolds},  arXiv:mathDG /0701723, 2007. 

\bibitem[5]{BBI} D. Burago, Y.  Burago, S. Ivanov, {\em A course in metric geometry}, 
Graduate Studies in Mathematics, 33. American Mathematical Society, Providence, RI, 2001. 


 \bibitem[6]{Ch} P. S. Chakraborty. {\em From C*-algebra extensions
     to CQMS, $SU_q(2), $ Podle{\`s} sphere and other examples. } 
arXiv:math/02100155v1.

\bibitem[7]{Co2} A. Connes. {\em Non Commutative Geometry}. Academic
    Press, San Diego, 1994.
    
\bibitem[8]{Co3}
A. Connes, {\em Gravity coupled with matter and the foundation of non
  commutative geometry},  Comm. Math. Phys.  {\bf 182}  (1996),  155--176.  

\bibitem[9]{CM}
A. Connes, H. Moscovici, {\em Transgression and the Chern character of finite-dimensional $K$-cycles}, Comm. Math. Phys. {\bf 155} (1993), 103--122.  


\bibitem[10]{DDLW}
L. D\c{a}browski, F.  D'Andrea, G. Landi, E. Wagner, {\em  Dirac operators on all Podle\'s quantum spheres},  J. Noncommut. Geom., {\bf 1} (2007), 213--239.  

\bibitem[11]{DS}
L. D\c{a}browski, A. Sitarz, {\em Dirac operator on the standard Podle\'s quantum sphere}, Noncommutative geometry and quantum groups (Warsaw, 2001), 49--58, Banach Center Publ., 61, Polish Acad. Sci., Warsaw, 2003.
    
  \bibitem[12]{Da}
K. R. Davidson, {\em $C^*$-algebras by example}, Fields Institute Monographs, {\bf 6}, American Mathematical Society, Providence, RI, 1996.


\bibitem[13]{Di} J. Dixmier. {\em Les C*-alg{\`e}bres et leurs
    repr{\'e}sentations}, Gauthier-Villars, Paris, 1964.

\bibitem[14]{Ef} E. G. Effros, {\em  Order ideals in a C*-algebra and
    its dual}, Duke Math. J., {\bf 30} (1963), 391--412. 
    
    
\bibitem[15]{Fu} K. Fukaya, {\em  Metric Riemannian Geometry}, 
Handbook of differential geometry, Vol. II, Elsevier/North-Holland, Amsterdam, 2006, 189--313. 


\bibitem[16]{HR} N.Higson, J. Roe, {\em Analytic K-Homology}, Oxford
  University Press, Oxford, 2000.

\bibitem[17]{Ke}  D. Kerr, {\em Matricial quantum Gromov-Hausdorff
    distance} J. Funct. Anal.,  {\bf 205}  (2003),  132--167.

\bibitem[18]{KL}  D. Kerr, H. Li {\em On Gromov-Hausdorff convergence
    for operator metric spaces}, arXiv: mathOA/0411157 v4, 2007. 
    
\bibitem[19]{KlL1} S. Klimek, A. Lesniewski, {\em A two-parameter quantum deformation of the unit disc}, 
 J. Funct. Anal., {\bf 115} (1993), 1--23.  


    
\bibitem[20]{KlL2} S. Klimek, A. Lesniewski, {\em Quantum Riemann surfaces. I. The unit disc}, 
 Comm. Math. Phys., {\bf 146} (1992), 103--122.  


\bibitem[21]{La} F. Latremoliere, {\em Approximation of quantum tori
    by finite quantum tori for the quantum Gromov-Hausdorff distance},
  J. Funct. Anal. {\bf 223 } (2005), 336--395.
 

 \bibitem[22]{Li}
H. Li, {\em Order-unit quantum Gromov-Hausdorff distance}, J. Funct. Anal. {\bf 231} (2006), 312--360.  

\bibitem[23]{Na} G. Nagy, {\em On the Haar measure of the quantum} ${\rm SU}(N)$ {\em group},  
Comm. Math. Phys., {\bf 153} (1993), 217--228. 


\bibitem[24]{Pe}
G.K. Pedersen, {\em Analysis now}, Graduate Texts in Mathematics, {\bf 118}, Springer-Verlag, New York, 1989.  


\bibitem[25]{Po}
P. Podle{\`s}, {\em Quantum spheres}, Lett. Math. Phys. {\bf
  14} 1987, 193--202.

\bibitem[26]{RiLsn}
 M. A. Rieffel,  {\em  Leibniz seminorms for "matrix algebras converge to the sphere''}, arXiv:mathOA/0707.3229, 2007.

  
 \bibitem [27]{RiMem}
 M. A. Rieffel, {\em Gromov--Hausdorff distance for quantum
   metric spaces},   Mem. Amer. Math. Soc.,  {\bf 168 }  (2004),
 no. 796, 1--65.   
 
 
\bibitem[28]{RiStSp}
M. A. Rieffel  {\em Metrics on state spaces},
Doc. Math. {\bf 4} (1999), 559--600.


\bibitem[29]{Sh} A. J. L. Sheu, {\em Quantization of the Poisson} ${\rm SU}(2)$ {\em and its Poisson homogeneous space---the} $2${\em -sphere. With an appendix by Jiang-Hua Lu and Alan Weinstein}, Comm. Math. Phys., {\bf  135} (1991), 217--232. 

\bibitem[30]{Wu} Wu. W, {\em Quantized Gromov-Hausdorff distance}, 
J. Funct. Anal. {\bf 238} (2006), 58--98.


\end{thebibliography}
\end{document}